\documentclass{amsart}

\newtheorem{theorem}[equation]{Theorem}
\newtheorem{lemma}[equation]{Lemma}
\newtheorem{corollary}[equation]{Corollary}
\newtheorem{proposition}[equation]{Proposition}

\numberwithin{equation}{section}

\usepackage{amscd}

\begin{document}

\title{On the Jacobian ring of a complete intersection}
\author{Alan Adolphson}
\address{Department of Mathematics\\
Oklahoma State University\\
Stillwater, Oklahoma 74078}
\email{adolphs@math.okstate.edu}
\author{Steven Sperber}
\address{School of Mathematics\\
University of Minnesota\\
Minneapolis, Minnesota 55455}
\email{sperber@math.umn.edu}
\date{\today}
\begin{abstract}
Let $f_1,\dots,f_r\in K[x]$, $K$ a field, be homogeneous polynomials and 
put $F=\sum_{i=1}^r y_if_i\in K[x,y]$.  The quotient $J=K[x,y]/I$, where
$I$ is the ideal generated by the $\partial F/\partial x_i$ and
$\partial F/\partial y_j$, is the {\it Jacobian ring of\/} $F$.  We
describe $J$ by computing the cohomology of a certain complex whose
top cohomology group is $J$.
\end{abstract}
\maketitle

\section{Introduction}

Let $K$ be a field and let $f_1,\dots,f_r\in K[x_1,\dots,x_n]$ be
homogeneous polynomials of degrees $d_1,\dots,d_r\geq 1$.  Set
\[ F=y_1f_1+\cdots+y_rf_r\in K[x_1,\dots,x_n,y_1,\dots,y_r]. \]
The quotient ring
\[ J=K[x,y]/(\partial F/\partial y_1,\dots,\partial F/\partial
y_r,\partial F/\partial x_1,\dots,\partial F/\partial x_n) \]
is sometimes referred to as the {\it Jacobian ring of\/} $F$.  (Note that
$\partial F/\partial y_j = f_j$.)  One reason for
interest in this ring is its connection with Hodge theory.  Consider
the bigrading $(\deg_1,\deg_2)$ on $K[x,y]$ defined by setting
\begin{align}
\deg_1 x_i &= 1,\;i=1,\dots,n, &\deg_2 x_i &= 0,\;i=1,\dots,n, \\
\deg_1 y_j&=-d_j,\;j=1,\dots,r, &\deg_2 y_j&=1,\;j=1,\dots,r.
\end{align}
Let $K[x,y]^{(q,p)}$ denote the bigraded component of bidegree
$(q,p)$.  Suppose that $K={\bf C}$, $r<n$, and the equations
$f_1=\cdots=f_r=0$ define a smooth complete intersection $X$ in
${\bf P}^{n-1}$ (i.~e., the $r\times n$ matrix with entries $\partial
f_i/\partial x_j$ has rank $r$ at every point of $X$).  Then
\begin{equation}
\dim_{\bf C}J^{(d_1+\cdots+d_r-n,p)} = \dim_{\bf C} H_{\rm
  prim}^{n-r-p-1}(X,\Omega_{X/{\bf C}}^p),
\end{equation}
where the subscript ``prim'' denotes the primitive subspace of the
cohomology group.  In this generality, this result is due to
Konno\cite{K} (see also Terasoma\cite{T}).  The hypersurface case is
due to Griffiths\cite{G}.  For further discussion, we refer the reader
to Dimca\cite{D}. 

The Jacobian ring also arises in the work of Dwork\cite{DW} and
Ireland\cite{I}.  Dwork showed that for smooth projective
hypersurfaces over a finite field $K$, a lower
bound for the Newton polygon of the interesting factor of the
zeta function is given by the polygon with sides of slopes
$p=0,1,\dots,n-r-1$ each with multiplicity $\dim_K J^{(d_1-n,p)}$.
Due to technical difficulties, Ireland was unable to completely extend
Dwork's result to the case of complete intersections.  Katz\cite{KA}
showed that the first side of the Newton polygon lies above the
conjectured lower bound.  The full result was later proved by
Mazur\cite{M}.  In \cite{AS1}, we used a toric approach to give
another proof; however, that approach required the stronger hypothesis
that the equation $f_1\cdots f_r=0$ define a normal crossing divisor
and that $d_1\cdots d_r\neq 0$ in $K$.

The purpose of this paper is to establish some results about $J$ that
will enable us in a future article to prove the theorem of Mazur by
generalizing Dwork's work on smooth projective hypersurfaces.  We
regard the Jacobian ring 
as the top cohomology group of a certain de Rham-type complex, namely,
the complex of differential forms $\Omega^{\bullet}_{K[x,y]/K}$ with
boundary map $\partial$ defined by $\partial(\omega)=dF\wedge\omega$.
Then clearly $J\cong H^{n+r}(\Omega^{\bullet}_{K[x,y]/K})$ as
$K[x,y]$-modules.  Each $\Omega^k_{K[x,y]/K}$ is given a bigrading by
extending the earlier bigrading: set 
\begin{align}
\deg_1 dx_i &= 1,\;i=1,\dots,n, &\deg_2 dx_i &= 0,\;i=1,\dots,n, \\
\deg_1 dy_j&=-d_j,\;j=1,\dots,r, &\deg_2 dy_j&=1,\;j=1,\dots,r.
\end{align}
Then $(\Omega^{\bullet}_{K[x,y]/K},\partial)$ is a bigraded complex
with boundary map of bidegree $(0,1)$.  In terms of this bigrading,
\[ J^{(d_1+\cdots+d_r-n,p)}\cong
H^{n+r}(\Omega^{\bullet}_{K[x,y]/K})^{(0,p+r)}. \]

The main problem we consider here is thus the computation of the
cohomology of the complexes
$(\Omega^{\bullet}_{K[x,y]/K},\partial)^{(0,\bullet)}$.  Note that in
our work we make no restriction on the characteristic of $K$, but
that an exceptional case arises when $d_1\cdots d_r=0$ in $K$.
Already in Dwork's treatment of hypersurfaces (\cite{DW1,DW}) a
similar exceptional case arose.  Our main result is the following.  
\begin{theorem}
Suppose that $r<n$ and that the equations $f_1=\cdots=f_r=0$ define a
smooth complete intersection $X$ in ${\bf P}^{n-1}$.  Then 
\begin{equation}
H^k(\Omega^{\bullet}_{K[x,y]/K})^{(0,p)} = 0\quad\text{for
  $k\neq 2r,n+r-1,n+r$ and all $p$} 
\end{equation}
and
\begin{equation}
H^{n+r}(\Omega^{\bullet}_{K[x,y]/K})^{(0,p)} = 0
\quad\text{if $p<r$ or $p\geq n$}.
\end{equation}
If $r<n-1$, then
\begin{equation}
H^{2r}(\Omega^{\bullet}_{K[x,y]/K})^{(0,p)} = \begin{cases}
K\cdot[df_1\wedge\cdots\wedge df_r\wedge dy_1\wedge\cdots\wedge dy_r]
& \text{if $p=r$} \\0 & \text{otherwise} \end{cases}
\end{equation}
and
\begin{equation}
\dim_K H^{n+r-1}(\Omega^{\bullet}_{K[x,y]/K})^{(0,p)} =
\dim_K H^{n+r}(\Omega^{\bullet}_{K[x,y]/K})^{(0,p)}
\end{equation}
unless $d_1\cdots d_r=0$ in $K$, $n+r$ is
even, and $p=\frac{n+r}{2}$ or $p=\frac{n+r}{2}-1$.  In this case,
\begin{multline}
\dim_K H^{n+r-1}(\Omega^{\bullet}_{K[x,y]/K})^{(0,p)} = \\
\dim_K H^{n+r}(\Omega^{\bullet}_{K[x,y]/K})^{(0,p)} +
\begin{cases} 1 & \text{if $p=\frac{n+r}{2}-1$} \\ -1 & \text{if
$p=\frac{n+r}{2}$.} \end{cases}
\end{multline}
If $r=n-1$ (so that $2r=n+r-1$), then
\begin{multline}
\dim_K H^{n+r-1}(\Omega^{\bullet}_{K[x,y]/K})^{(0,p)} = \\
\dim_K H^{n+r}(\Omega^{\bullet}_{K[x,y]/K})^{(0,p)} + \begin{cases} 1 &
\text{if $p=r$} \\ 0 & \text{otherwise.} \end{cases}
\end{multline}
\end{theorem}

Put
\begin{multline*}
h_p = \dim_K H^{n+r}(\Omega^{\bullet}_{K[x,y]/K})^{(0,p)} + \\
\begin{cases} -1 & \text{if $d_1\cdots d_r=0$ in $K$, $n+r$ 
is even, and $p=\frac{n+r}{2}$} \\ 0 & \text{otherwise,} \end{cases} 
\end{multline*}
and define a polynomial $H(t)$ by
\begin{equation}
H(t) = \sum_{p=r}^{n-1} h_pt^p.
\end{equation}
Then $H(t)$ (or $H(t) + t^{(n+r)/2}$, in the exceptional case) is the
Hilbert series of the graded module
$H^{n+r}(\Omega^{\bullet}_{K[x,y]/K})^{(0,\bullet)}$ (see
\cite[Section 13]{MA} for general information on Hilbert series).
Since $\dim_K (\Omega^k_{K[x,y]/K})^{(0,p)}$ 
is easily expressible as a binomial coefficient, one obtains a formula
for $H(t)$ from Theorem~1.6 (see equation (5.20) below) that shows
$H(t)$ is independent of $K$.  From this formula it is straightforward  
to calculate $H(1)$ and to check that $t^{n+r-1}H(1/t)=H(t)$, giving
the usual formula for the dimension of the primitive part of
middle-dimensional cohomology and proving the symmetry of the Hodge
numbers.  This computation will be carried out in section 5 to give
the following result.
\begin{corollary}
Under the hypotheses of{\/} {\rm Theorem}~$1.6$, 
\begin{equation} 
\sum_{p=r}^{n-1}h_p = (-1)^{n-r}(n-r) + (-1)^n \sum_{l=r}^{n-1} (-1)^{l+1}
\binom{n}{l+1} \sum_{\substack{i_1+\cdots+i_r=l \\ i_j\geq 1 \text{for
all $j$}}} d_1^{i_1}\cdots d_r^{i_r} 
\end{equation}
and
\begin{equation}
h_p = h_{n+r-1-p}\quad\text{for all $p$.}
\end{equation}
\end{corollary}

A key technical step is Proposition~2.2, which is a complement to the
well known de Rham-Saito Lemma.  The de Rham-Saito Lemma and
Proposition~2.2 are special cases of a more general result,
Theorem~2.15, which is not used in this paper but which is included
for completeness. 
The complexes $(\Omega^\bullet_{K[x,y]/K},\partial)^{(0,\bullet)}$ are
defined even when $r\geq n$.  In future work, we plan to relate them
to the complement of the divisor $f_1\cdots f_r=0$ in ${\bf P}^{n-1}$.
In the case where the ideal $(f_1,\dots,f_r)$ has depth $n$, we sketch
the computation of the cohomology of these complexes in section 6.   

To facilitate the calculation of
$H^{n+r-1}(\Omega^\bullet_{K[x,y]/K})^{(0,\bullet)}$ in section 4, we 
introduce a subcomplex
$(\widetilde{\Omega}^\bullet_{K[x,y]/K},\partial)^{(0,\bullet)}$ of
$(\Omega^\bullet_{K[x,y]/K},\partial)^{(0,\bullet)}$ and compute its
cohomology.  From this computation one sees that when $K={\bf C}$
\[ \bigoplus_p H^i(\widetilde{\Omega}^\bullet_{{\bf C}[x,y]/{\bf
C}})^{(0,p)} \cong H^{i-1}_{\rm DR}({\bf P}^{n-1}\setminus X). \] 
The complex $(\widetilde{\Omega}^\bullet_{{\bf C}[x,y]/{\bf
C}},d+dF\wedge)$, where the boundary operator $\partial$ has been
replaced by $d+dF\wedge$, is no longer bigraded but remains graded
relative to $\deg_1$.  In a future article, we shall show that for
arbitrary homogeneous polynomials $f_1,\dots,f_r$ there is a
quasi-isomorphism from $(\widetilde{\Omega}^\bullet_{{\bf C}[x,y]/{\bf
C}},d+dF\wedge)^{(0)}$ to the usual \v{C}ech-de Rham
complex of ${\bf P}^{n-1}\setminus X$ relative to the collection of
open sets defined by $f_i\neq 0$, $i=1,\dots,r$.  Furthermore, when
$X$ is a smooth complete intersection, the filtration $\deg_2$
on $(\widetilde{\Omega}^\bullet_{{\bf C}[x,y]/{\bf
C}},d+dF\wedge)^{(0)}$ is identified to the Hodge filtration on
this \v{C}ech-de Rham complex under this quasi-isomorphism.  Passing
to the associated graded complexes leads to another proof of the
relation (1.3).  

We are indebted to G. Lyubeznik and C. Huneke for pointing out the
reference \cite{SS} in connection with section 6.

\section{A complement to the de Rham-Saito Lemma}

We begin by reminding the reader of the de Rham-Saito Lemma\cite{S}.  Let
$A$ be a commutative Noetherian ring with $1$ and let $M$ be a free
$A$-module of rank $n$ with basis $e_1,\ldots,e_n$.  We denote
its $k$-th exterior power by $\wedge^k M$.  This is a free
$A$-module with basis 
\[ \{e_{i_1}\wedge\cdots\wedge e_{i_k} \mid 1\leq i_1<\cdots<i_k\leq
n\}. \]
Fix $\omega_1,\ldots,\omega_r\in M$ and write
\[ \omega_1\wedge\cdots\wedge\omega_r=\sum_{1\leq i_1<\cdots<i_r\leq
n} a_{i_1\cdots i_r}\,e_{i_1}\wedge\cdots\wedge e_{i_r}. \]
Let $I$ be the ideal of $A$ generated by the $a_{i_1\cdots i_r}$.  We
note for future use that for every subset $\{i_1,\ldots,i_k\}$ of
$\{1,\ldots,r\}$, the ideal generated by the coefficients of
$\omega_{i_1}\wedge\cdots\wedge\omega_{i_k}$ contains $I$.  The
main result of \cite{S} is the following.
\begin{proposition}
Let $\omega\in \wedge^k M$ satisfy
$\omega_1\wedge\cdots\wedge\omega_r\wedge\omega=0$.  \\
{\bf (a)} There exists $m\geq 0$ with the property that if $g\in I^m$, then
\[ g\omega=\sum_{i=1}^r \omega_i\wedge\alpha_i \]
for some $\alpha_1,\ldots,\alpha_r\in\wedge^{k-1}M$. \\
{\bf (b)} If $k<{\rm depth}(I)$, then there exist
$\alpha_1,\ldots,\alpha_r\in\wedge^{k-1}M$ such that 
\[ \omega=\sum_{i=1}^r \omega_i\wedge\alpha_i. \]
\end{proposition}

To analyze the Jacobian ring, we begin with the following result.
\begin{proposition}
Let $\omega\in \wedge^k M$ satisfy $\omega_i\wedge\omega=0$
for $i=1,\ldots,r$.  \\
{\bf (a)} There exists $m\geq 0$ with the property that if $g\in I^m$, then
\[ g\omega=\omega_1\wedge\cdots\wedge\omega_r\wedge\alpha \]
for some $\alpha\in\wedge^{k-r}M$.  \\
{\bf (b)} If ${\rm depth}(I)>0$ and $k<{\rm depth}(I)+r-1$, then there
exists $\alpha\in\wedge^{k-r}M$ such that 
\[ \omega=\omega_1\wedge\cdots\wedge\omega_r\wedge\alpha. \]
\end{proposition}

{\it Proof}.  Fix $g\in I$.  Since $\omega_1\wedge\omega=0$,
Proposition 2.1(a) implies there exists $m_1\geq 0$ and $\alpha_1\in
\wedge^{k-1}M$ such that $g^{m_1}\omega=\omega_1\wedge\alpha_1$.
Proceeding inductively, suppose that for some $i$, $1\leq i<r$, we have
\begin{equation}
g^m\omega=\omega_1\wedge\cdots\wedge\omega_i\wedge\beta 
\end{equation}
for some $m\geq 0$ and $\beta\in\wedge^{k-i}M$.  Since
$\omega_{i+1}\wedge\omega=0$, we get
\[ \omega_1\wedge\cdots\wedge\omega_i\wedge\omega_{i+1}\wedge\beta
=0, \]
so by Proposition 2.1(a) there exist $m_i\geq 0$ and
$\beta_1,\ldots,\beta_{i+1}\in \wedge^{k-i-1} M$ such that
\[ g^{m_i}\beta=\sum_{j=1}^{i+1} \omega_j\wedge\beta_j. \]
Substitution into (2.3) then gives
\[
g^{m+m_i}\omega=\omega_1\wedge\cdots\wedge\omega_{i+1}\wedge\beta_{i+1}.
\]
By induction, we arrive at the equation
\begin{equation}
g^{m'}\omega=\omega_1\wedge\cdots\wedge\omega_r\wedge\gamma
\end{equation}
for some $m'\geq 0$ and some $\gamma\in\wedge^{k-r}M$.  Now $m'$ and
$\gamma$ depend on $g$, but $I$ is finitely generated, say, by
$g_1,\dots,g_s$.  It follows that there exists an integer $m(I)$ and
$\gamma_1,\dots,\gamma_s\in\wedge^{k-r}M$ such that
\[ g_i^{m(I)}\omega = \omega_1\wedge\cdots\wedge\omega_r\wedge
\gamma_i. \]
Since $I^{m(I)-s+1}\subseteq(g_1^{m(I)},\dots,g_s^{m(I)})$, this
establishes part (a) of the proposition.  

For part (b), again fix $g\in I$ and note that (2.4) implies 
\begin{equation}
\omega_1\wedge\cdots\wedge\omega_r\wedge\gamma=0 \quad\text{in
$\wedge^k M/g^{m'}M$}.
\end{equation}
Since we are assuming $\text{depth}(I)>0$, we may assume $g$ is not a
zero-divisor in~$A$.  The ideal $I/(g^{m'})$ then has depth equal to
$\text{depth}(I)-1$.  Proposition 2.1(b) now says that for
$k-r<\text{depth}(I)-1$ (i.~e., for $k<\text{depth}(I)+r-1$), there
exist $\gamma_1,\ldots,\gamma_r\in \wedge^{k-r-1}M$ and
$\gamma_0\in\wedge^{k-r}M$ such that 
\[ \gamma=\sum_{i=1}^r \omega_i\wedge\gamma_i + g^{m'}\gamma_0. \]
Substituting this expression into (2.4) gives
\[ g^{m'}(\omega-\omega_1\wedge\cdots\wedge\omega_r\wedge\gamma_0) =
0, \]
and since $g$ is not a zero-divisor in $A$ we conclude that
\[ \omega=\omega_1\wedge\cdots\wedge\omega_r\wedge\gamma_0, \]
which establishes part (b) of the proposition.

We apply Propositions 2.1 and 2.2 to prove a result that is more
directly connected with the Jacobian ring.  Let $y_1,\ldots,y_r$ be 
indeterminates and consider
\[ M'=A[y_1,\ldots,y_r]\otimes_A M, \]
a free $A[y_1,\ldots,y_r]$-module.  For simplicity, we write its basis
as $e_1,\ldots,e_n$ instead of $1\otimes e_1,\ldots,1\otimes e_n$ and
we write $\omega_1,\ldots,\omega_r$ instead of $1\otimes\omega_1,
\ldots,1\otimes\omega_r$.  For all~$k$,
\[ \wedge^k M'=A[y_1,\ldots,y_r]\otimes_A (\wedge^k M) \]
and, as $A$-modules,
\[ \wedge^k M'= \bigoplus_{b_1,\ldots,b_r\geq 0} y_1^{b_1}\cdots
y_r^{b_r}(\wedge^k M). \]
There is thus a natural grading on $\wedge^k M'$ by $\text{deg}_y$,
the total degree in $y_1,\ldots,y_r$.  We denote by $(\wedge^k
M')^{(p)}$ the homogeneous component of degree $p$ in this grading,
\[ (\wedge^k M')^{(p)} = \bigoplus_{b_1+\cdots+b_r=p} y_1^{b_1}\cdots
y_r^{b_r}(\wedge^k M), \]
and we make the identification $(\wedge^k M')^{(0)}=\wedge^k M$.  

Consider the element $\sum_{i=1}^r y_i\omega_i\in (\wedge^1
M')^{(1)}$ and define $\partial:\wedge^k M'\rightarrow \wedge^{k+1} M'$ by
\[ \partial(\omega)=(y_1\omega_1+\cdots+y_r\omega_r)\wedge\omega. \]
Note that $\partial$ is homogeneous of degree 1, hence this is a
graded complex.  We denote by $H^k(\wedge^\bullet M',\partial)^{(p)}$
the homogeneous component of degree $p$ in the induced grading on
cohomology. 

If $\omega\in (\wedge^k M')^{(p)}$, we may write
\begin{equation}
\omega=\sum_{b_1+\cdots+b_r=p} y_1^{b_1}\cdots y_r^{b_r}
\omega(b_1,\ldots,b_r)
\end{equation}
with $\omega(b_1,\ldots,b_r)\in \wedge^k M$.  The condition
$\partial(\omega)=0$ is equivalent to the equations
\begin{equation}
\sum_{j=1}^r
\omega_j\wedge\omega(c_1,\ldots,c_{j-1},c_j-1,c_{j+1},
\ldots,c_r)=0
\end{equation}
for all nonnegative integers $c_1,\ldots,c_r$ satisfying
$c_1+\cdots+c_r=p+1$ (with the understanding that
$\omega(b_1,\ldots,b_r)=0$ if $b_i<0$ for some $i$).  Let $\alpha\in
(\wedge^{k-1} M')^{(p-1)}$, say,
\begin{equation}
\alpha=\sum_{a_1+\cdots+a_r=p-1} y_1^{a_1}\cdots y_r^{a_r}
\alpha(a_1,\ldots,a_r) 
\end{equation}
with $\alpha(a_1,\ldots,a_r)\in\wedge^{k-1}M$.  The condition
$\partial(\alpha)=\omega$ is equivalent to the equations
\begin{equation}
\omega(b_1,\ldots,b_r)=\sum_{j=1}^r \omega_j\wedge
\alpha(b_1,\ldots,b_{j-1}, b_j-1,b_{j+1},\ldots,b_r)
\end{equation}
for all nonnegative integers $b_1,\ldots,b_r$ satisfying
$b_1+\cdots+b_r=p$.  

\begin{proposition} 
{\bf (a)} For each $p\geq 1$ and $k\geq 0$ there exists $m=m_{k,p}\geq
0$ such that 
\[ I^m H^k({\wedge}^\bullet M',\partial)^{(p)} = 0. \]
{\bf (b)} If $1\leq p<{\rm depth}(I)$ and $k<{\rm depth}(I)+r-1$, then
\[ H^k({\wedge}^\bullet M',\partial)^{(p)} = 0. \]
\end{proposition}

{\it Remark}.  Note that Proposition 2.2 may be regarded as describing
this cohomology when $p=0$.  Proposition 2.2(a) is equivalent to the
assertion that
\[ I^m H^k({\wedge}^\bullet M',\partial)^{(0)}
\subseteq \omega_1\wedge\cdots\wedge\omega_r\wedge(\wedge^{k-r} M) \]
for $m$ sufficiently large, while Proposition 2.2(b) is equivalent to
the assertion that, if $\text{depth}(I)>0$ and
$k<\text{depth}(I)+r-1$, then
\[ H^k({\wedge}^\bullet M',\partial)^{(0)} = 
\omega_1\wedge\cdots\wedge\omega_r\wedge(\wedge^{k-r} M). \]

{\it Proof}.  We prove part (a) by induction on $p$.  We establish the
case $p=1$ by induction on $r$.  For $r=1$ the assertion
follows immediately from Proposition 2.1(a), so assume the result true for
$r-1$.  Let $\omega\in (\wedge^kM')^{(1)}$ with
$\partial(\omega)=0$.  Let $\epsilon_i$ be the $r$-tuple with 1 in the
$i$-th position and zeros elsewhere.  By (2.6) we may write
\[ \omega=y_1\omega(\epsilon_1)+\cdots+y_r\omega(\epsilon_r) \]
with $\omega(\epsilon_i)\in\wedge^k M$.  Since 
\[ (y_1\omega_1+\cdots+y_r\omega_r)\wedge
(y_1\omega(\epsilon_1)+\cdots+y_r\omega(\epsilon_r)) = 0, \]
we also have
\[ (y_1\omega_1+\cdots+y_{r-1}\omega_{r-1})\wedge
(y_1\omega(\epsilon_1)+\cdots+y_{r-1}\omega(\epsilon_{r-1})) = 0. \]
By induction, there exists $m\geq 0$ such that for all $g\in I^m$, 
\[ g(y_1\omega(\epsilon_1)+\cdots+y_{r-1}\omega(\epsilon_{r-1})) = 
(y_1\omega_1+\cdots+y_{r-1}\omega_{r-1})\wedge \alpha \]
for some $\alpha\in\wedge^{k-1}M$.  It follows that
\begin{equation}
g\omega-\partial(\alpha)=y_r(g\omega(\epsilon_r)-\omega_r\wedge
\alpha), 
\end{equation}
so the equation $\partial(g\omega-\partial(\alpha))=0$ reduces to
\[ \omega_i\wedge(g\omega(\epsilon_r)-\omega_r\wedge\alpha)=0
\qquad\text{for $i=1,\ldots,r$.} \]
We apply Proposition 2.2(a) to conclude there exists $m'\geq 0$ such
that for $h\in I^{m'}$,
\[ h(g\omega(\epsilon_r)-\omega_r\wedge\alpha)=\omega_1\wedge\cdots
\wedge\omega_r\wedge\beta \]
for some $\beta\in\wedge^{k-r}M$.  Equation (2.11) then implies
\[ (gh)\omega=\partial(h\alpha+(-1)^{r-1}\omega_1\wedge\cdots\wedge
\omega_{r-1}\wedge\beta), \]
i.~e., $I^{m+m'}$ annihilates $H^k(\wedge^\bullet
M',\partial)^{(1)}$.  

Now suppose the assertion true for $p-1$.  We prove it for $p$ by
induction on $r$.  For $r=1$, the assertion follows immediately from
Proposition 2.1(a), so assume the result true for $r-1$.  Let
$\omega\in(\wedge^k M')^{(p)}$ with $\partial(\omega)=0$ and let
$\omega$ be written as in~(2.6).  Put
\[ \omega'=\sum_{b_1+\cdots+b_{r-1}=p} y_1^{b_1}\cdots
y_{r-1}^{b_{r-1}} \omega(b_1,\ldots,b_{r-1},0). \]
The condition $\partial(\omega)=0$ implies
\[ (y_1\omega_1+\cdots+y_{r-1}\omega_{r-1})\wedge\omega'=0, \]
so by induction on $r$ there exists $m\geq 0$ such that if $g\in I^m$,
then there exists 
\[ \alpha=\sum_{a_1+\cdots+a_{r-1}=p-1} y_1^{a_1}\cdots
y_{r-1}^{a_{r-1}} \alpha(a_1,\ldots,a_{r-1}), \]
with $\alpha(a_1,\ldots,a_{r-1})\in\wedge^{k-1}M$, such that
\[ g\omega'=(y_1\omega_1+\cdots+y_{r-1}\omega_{r-1})\wedge\alpha. \]
It follows that all terms of $g\omega-\partial(\alpha)$ are divisible
by $y_r$, say,
\begin{equation}
g\omega-\partial(\alpha)=y_r\beta
\end{equation}
for some $\beta\in (\wedge^kM')^{(p-1)}$.  This equation implies
$\partial(\beta)=0$, so by induction on $p$ there exists $m'\geq 0$
such that if $h\in I^{m'}$, then
\[ h\beta=\partial(\gamma) \]
for some $\gamma\in(\wedge^{k-1}M')^{(p-2)}$.  Equation (2.12) then
implies
\[ (gh)\omega=\partial(h\alpha+y_r\gamma), \]
i.~e., $I^{m+m'}$ annihilates $H^k(\wedge^\bullet
M',\partial)^{(p)}$.  This completes the proof of part (a).

To prove part (b), we begin with some notation.  In the
course of the argument, we shall produce a sequence $g_1,\ldots,g_p\in
I$.  For $i=1,\ldots,p$, put
\begin{align*}
A_i &=A/(g_1,\ldots,g_i), \\
M_i &=M/(g_1,\ldots,g_i)M  \\
M'_i &=M'/(g_1,\ldots,g_i)M'.
\end{align*}
These are free $A_i$-modules.  The boundary map $\partial:\wedge^k
M'\rightarrow \wedge^{k+1}M'$ induces $\partial:\wedge^k
M'_i\rightarrow \wedge^{k+1}M'_i$ for $i=1,\ldots,p$, producing
complexes $(\wedge^\bullet M'_i,\partial)$.  

Fix $k$ and~$p$ satisfying the hypothesis of part (b) and let
$\alpha_0\in (\wedge^k M')^{(p)}$ with $\partial(\alpha_0)=0$.  By
part (a), for all $g_1$ in some power of $I$ there exists $\alpha_1\in
(\wedge^{k-1}M')^{(p-1)}$ such that
\[ g_1\alpha_0=\partial(\alpha_1). \]
Suppose inductively
that for some $i$, $1\leq i<p$, we have chosen $g_i$ in some power 
of~$I$ and $\alpha_i\in (\wedge^{k-i}M'_{i-1})^{(p-i)}$ such that 
\begin{equation}
g_i\alpha_{i-1}=\partial(\alpha_i)\qquad\text{in $\wedge^{k-i+1}
M'_{i-1}$.}
\end{equation}
Then by part (a) applied to the complex $(\wedge^\bullet
M'_i,\partial)$, for all $g_{i+1}$ in some power of~$I$ there 
exists $\alpha_{i+1}\in (\wedge^{k-i-1}M'_i)^{(p-i-1)}$ such that
\[ g_{i+1}\alpha_i=\partial(\alpha_{i+1})\qquad\text{in
$\wedge^{k-i} M'_i$.} \]
Thus by induction equation (2.13) holds for $i=1,\ldots,p$.  

Note that since $p<\text{depth}(I)$, we may choose $g_1,\ldots,g_p$ to
be a regular sequence in $A$.  Taking $i=p$ in (2.13) gives
\[ \alpha_p\in (\wedge^{k-p}M_{p-1}')^{(0)}=\wedge^{k-p}M_{p-1} \] 
satisfying
\[ \partial(\alpha_p)=0\qquad\text{in $\wedge^{k-p+1}M'_p$.} \]
This is equivalent to the condition that $\omega_i\wedge\alpha_p=0$ in
$\wedge^{k-p+1}M_p$ for $i=1,\ldots,r$.  Since
$\text{depth}(I/(g_1,\ldots,g_p))=\text{depth}(I)-p>0$ and $k-p<{\rm
depth}(I/(g_1,\ldots,g_p)) + r-1$, we may
apply Proposition 2.2(b) in the ring $A_p$ to conclude that 
\[ \alpha_p=\omega_1\wedge\cdots\wedge\omega_r\wedge\beta_p +
g_p\gamma_p \]
for some $\beta_p\in \wedge^{k-p-r}M_{p-1}$, $\gamma_p\in
\wedge^{k-p}M_{p-1}$.  If we now take $i=p$ in (2.13) and
substitute this expression for $\alpha_p$ we get
\[ g_p\alpha_{p-1}=\partial(g_p\gamma_p) \qquad\text{in
$\wedge^{k-p+1}M'_{p-1}$.} \]
Since $g_p$ is not a zero-divisor in $A_{p-1}$, this implies
\[ \alpha_{p-1}=\partial(\gamma_p)\qquad\text{in
$\wedge^{k-p+1}M'_{p-1}$.} \] 
Suppose that for some $i$, $1\leq i\leq p-1$, we have shown
\begin{equation}
\alpha_i=\partial(\gamma_{i+1})\qquad\text{in $\wedge^{k-i}M'_i$}
\end{equation}
for some $\gamma_{i+1}\in\wedge^{k-i-1}M'_i$.  Let
$\tilde{\gamma}_{i+1}$ be any lifting of $\gamma_{i+1}$ to
$\wedge^{k-i-1}M'_{i-1}$.  Then (2.14) implies
\[ \alpha_i=\partial(\tilde{\gamma}_{i+1}) + g_i\gamma_i
\qquad\text{in $\wedge^{k-i}M'_{i-1}$} \]
for some $\gamma_i\in\wedge^{k-i}M'_{i-1}$.  Substitution into
(2.13) gives
\[ g_i\alpha_{i-1}=\partial(g_i\gamma_i)\qquad\text{in
$\wedge^{k-i+1}M'_{i-1}$.} \]
Since $g_i$ is not a zero-divisor in $A_{i-1}$, we conclude that
(2.14) holds with $i$ replaced by $i-1$.  By descending induction on
$i$, it follows that (2.14) holds for $i=0$, which is the assertion of
part (b).

We give a generalization of Propositions 2.1 and 2.2.  This result is
included for the sake of completeness and is not used in this article.
Let $N$ denote the $A$-submodule of $M$ spanned by
$\omega_1,\dots,\omega_r$.  For
$J=\{j_1,\dots,j_t\}\subseteq\{1,\dots,r\}$, $j_1<\dots<j_t$, set
$\omega_J= \omega_{j_1}\wedge\cdots\wedge \omega_{j_t}$.  
\begin{theorem}
Fix $s$, $1\leq s\leq r$, and let $\omega\in\wedge^k M$ satisfy
$\gamma\wedge\omega=0$ for all $\gamma\in\wedge^s N$. \\
{\bf (a)} There exists $m\geq 0$ with the property that if $g\in I^m$, then
\[ g\omega = \sum_{\substack{J\subseteq\{1,\dots,r\} \\ |J|=r-s+1}}
\omega_J\wedge\alpha_J \]
for some $\alpha_J\in \wedge^{k-r+s-1}M$.  \\
{\bf (b)} If ${\rm depth}(I)>0$ and $k<{\rm depth}(I)+r-s$, then there exist
$\alpha_J\in \wedge^{k-r+s-1}M$ such that
\[ \omega = \sum_{\substack{J\subseteq\{1,\dots,r\} \\ |J|=r-s+1}}
\omega_J\wedge\alpha_J \]
\end{theorem}

{\it Proof}.  Since
\[ \omega_1\wedge\cdots\wedge\omega_s\wedge\omega = 0 \]
we have by Proposition 2.1(a) that there exists $m_1\geq 0$ such that
if $g\in I$, then there exist $\alpha_i\in\wedge^{k-1}M$ for
$i=1,\dots,s$ such that
\[ g^{m_1}\omega = \sum_{i=1}^s \omega_i\wedge\alpha_i. \]
Suppose that for some $t$, $1\leq t<r-s+1$, we have shown that there
exists $m_t\geq 0$ such that if $g\in I$, then 
\begin{equation}
g^{m_t}\omega = \sum_{\substack{B\subseteq\{1,\dots,s+t-1\} \\ |B|=t}}
\omega_B\wedge\alpha_B 
\end{equation}
for some $\alpha_B\in\wedge^{k-t}M$.  To complete the proof of part
(a) of the theorem, it suffices by induction on $t$ to prove that this
equation holds with $t$ replaced by $t+1$.  Fix
$C\subseteq\{1,\dots,s+t-1\}$ with $|C|=t$ and let $C'$ be the
complement of $C$ in $\{1,\dots,s+t\}$.  Take the wedge product
of both sides of (2.16) with $\omega_{C'}\in \wedge^s N$.  
By the hypothesis on $\omega$ we get
\begin{equation}
\omega_{C'}\wedge\sum_{\substack{B\subseteq\{1,\dots,s+t-1\} \\
|B|=t}} \omega_B\wedge\alpha_B = 0.  
\end{equation}
Note that the sets $C'$ and $B$ are not disjoint unless $B=C$, hence
$\omega_{C'}\wedge\omega_B = 0$ unless $B=C$.  It then follows from
(2.17) that 
\[ \omega_1\wedge\cdots\wedge\omega_{s+t}\wedge\alpha_B=0 \]
for every $B\subseteq\{1,\dots,s+t-1\}$, $|B|=t$.  Applying
Proposition 2.1(a), we conclude that there exists $m_B\geq 0$ such
that if $g\in I$, then
\begin{equation}
g^{m_B}\alpha_B = \sum_{i=1}^{s+t}\omega_i\wedge\alpha_{B,i}.
\end{equation}
Substituting (2.18) into (2.16) gives (2.16) with $t$ replaced by
$t+1$.

We prove part (b) by induction on $s$.  The case $s=1$ is Proposition
2.2(b).  So we assume the result holds for some $s$, $1\leq s<r$, and
prove it for $s+1$.  Thus we assume $\omega\in\wedge^k M$, $k<{\rm
depth}(I)+r-s-1$, and $\gamma\wedge\omega=0$ for all
$\gamma\in\wedge^{s+1}N$.  By part (a) of the theorem, we know there
exists $m\geq 0$ such that if $g\in I^m$, then
\begin{equation}
g\omega = \sum_{\substack{J\subseteq\{1,\dots,r\} \\ |J|=r-s}}
\omega_J\wedge\alpha_J
\end{equation}
for some $\alpha_J\in \wedge^{k-r+s}M$.  Since ${\rm depth}(I)>0$ we
may assume that $g$ is not a zero-divisor in $A$.  Equation (2.19)
implies 
\begin{equation}
\sum_{\substack{J\subseteq\{1,\dots,r\} \\ |J|=r-s}}
\omega_J\wedge\alpha_J = 0 \quad\text{in $\wedge^k
(M/gM)$.}
\end{equation}
Fix $B\subseteq\{1,\dots,r\}$, $|B|=r-s$, and let $B'$ be the
complement of $B$ in $\{1,\dots,r\}$.  Take the wedge product of
both sides of (2.20) with $\omega_{B'}$.  For
$J\subseteq\{1,\dots,r\}$, $|J|=r-s$, the sets $B'$ and $J$ are disjoint if and only if $J=B$,
hence (2.20) implies
\[ \omega_1\wedge\cdots\wedge\omega_r\wedge \alpha_B=0\quad\text{in $\wedge^k
(M/gM)$.} \]
Now $\alpha_B$ is a ($k-r+s$)-form and our hypothesis implies
\[ k-r+s<{\rm depth}(I)-1={\rm depth}(I/(g)) \]
so we may apply Proposition 2.1(b) to get
\[ \alpha_B=\sum_{i=1}^r \omega_i\wedge\alpha_{B,i} + g\beta_B, \]
where $\alpha_{B,i}\in\wedge^{k-r+s-1}M$ and
$\beta_B\in\wedge^{k-r+s}M$.  Substituting into (2.19) gives
\begin{equation}
g\bigg(\omega - \sum_{\substack{J\subseteq\{1,\dots,r\} \\ |J|=r-s}}
\omega_J\wedge\beta_J\bigg) = 
\sum_{\substack{J\subseteq\{1,\dots,r\} \\ |J|=r-s}} \omega_J\wedge
\sum_{i=1}^r \omega_i \wedge\alpha_{J,i}. 
\end{equation}
Let $C\subseteq\{1,\dots,r\}$, $|C|=s$.  Then
$\omega_C\wedge\omega_J\wedge\omega_i = 0$ for all
$J\subseteq\{1,\dots,r\}$, $|J|=r-s$, and all $i=1,\dots,r$, 
since $\wedge^{r+1}N =0$.  It follows that taking the wedge product
with $\omega_C$ annihilates the left-hand side of (2.21).  Since $g$
is not a zero divisor in $A$, we conclude that 
\[ \omega_C\wedge\bigg(\omega -
\sum_{\substack{J\subseteq\{1,\dots,r\} \\ |J|=r-s}}
\omega_J\wedge\beta_J\bigg) = 0 \] 
for all $C\subseteq\{1,\dots,r\}$, $|C|=s$.  Since
\[ k<{\rm depth}(I)+r-s-1<{\rm depth}(I)+r-s, \]
we may apply the induction hypothesis on $s$ to conclude that
\[ \omega - \sum_{\substack{J\subseteq\{1,\dots,r\} \\ |J|=r-s}}
\omega_J\wedge\beta_J =
\sum_{\substack{B\subseteq\{1,\dots,r\} \\ |B|=r-s+1}} 
\omega_B\wedge\beta_B \]
for some $\beta_B\in\wedge^{k-r+s-1}M$.  Solving this equation for
$\omega$ gives the assertion for $s+1$.

\section{Complete intersections and the Jacobian ring}

In this section, we prove assertions (1.7), (1.8), and (1.9) of
Theorem~1.6.  Returning to the situation described in the
introduction, let $K$ be
a field and let $f_1,\ldots,f_r\in K[x_1,\ldots,x_n]$, $1\leq r< n$,
be homogeneous polynomials of degrees $d_1,\ldots,d_r\geq 1$.  Set
\[ F=y_1f_1+\cdots+y_rf_r\in K[x_1,\ldots,x_n,y_1,\ldots,y_r]. \]
We are interested in studying the cohomology of the complex of
$K[x,y]$-modules $\Omega^\bullet_{K[x,y]/K}$ with boundary operator
$\partial:\Omega^k_{K[x,y]/K}\to\Omega^{k+1}_{K[x,y]/K}$ defined by
$\partial(\omega) = dF\wedge\omega$, where $dF\in \Omega^1_{K[x,y]/K}$
is the exterior derivative of $F$.  This complex is bigraded by the
bigrading defined in (1.1), (1.2), (1.4), and (1.5), and the boundary map
$\partial$ has bidegree $(0,1)$.

It is convenient to regard $(\Omega^\bullet_{K[x,y]/K},\partial)$ as
the total complex associated to the double complex
$(C^{l,m},\partial_h,\partial_v)$, where 
\[ C^{l,m} = \bigoplus_{\substack{1\leq i_1<\cdots<i_l\leq n \\ 1\leq
    j_1<\cdots<j_m\leq r}} K[x,y]\,dx_{i_1}\wedge\cdots\wedge
    dx_{i_l}\wedge dy_{j_1}\wedge\cdots\wedge dy_{j_m} \]
and $\partial_h:C^{l,m}\to C^{l+1,m}$ and $\partial_v:C^{l,m}\to
    C^{l,m+1}$ are defined by
\[ \partial_h(\omega) = (y_1df_1+\cdots+y_rdf_r)\wedge\omega,
    \qquad \partial_v(\omega) = (f_1dy_1+\cdots+f_rdy_r)\wedge\omega. \]

When $f_1,\ldots,f_r$ form a regular sequence in $K[x]$, the
cohomology of each column $(C^{l,\bullet},\partial_v)$ vanishes except in
dimension $r$, where one has
\[ H^r(C^{l,\bullet},\partial_v) = \bigoplus_{1\leq i_1<\cdots<i_l\leq n}
K[x,y]/(f_1,\ldots,f_r)\,dx_{i_1}\wedge\cdots\wedge dx_{i_l}\wedge
dy_1\wedge\cdots\wedge dy_r. \]
If we set $A=K[x]/(f_1,\ldots,f_r)$,
$M=A\otimes_{K[x]}\Omega^1_{K[x]/K}$, and
\[ M'=A[y_1,\ldots,y_r]\otimes_{K[x]}\Omega^1_{K[x]/K}, \]
we can write this more compactly as
\begin{equation}
H^r(C^{l,\bullet},\partial_v) = (\wedge^l M')\wedge dy_1\wedge\cdots\wedge
dy_r.
\end{equation}
It follows by a well-known result in commutative algebra that
\[ H^{k+r}(\Omega^\bullet_{K[x,y]/K}) = H^k((\wedge^\bullet
M')\wedge dy_1\wedge\cdots\wedge dy_r,\bar{\partial}_h), \]
where $\bar{\partial}_h:(\wedge^kM')\wedge dy_1\wedge\cdots\wedge
dy_r\to (\wedge^{k+1}M') \wedge dy_1\wedge\cdots\wedge dy_r$ is the map 
induced by $\partial_h$.  In particular, we conclude that
\begin{equation}
H^k(\Omega^\bullet_{K[x,y]/K}) = 0 \quad\text{for $k<r$.}
\end{equation}

It is notationally convenient to drop the symbol
``$dy_1\wedge\cdots\wedge dy_r$'' and adjust the bigrading
accordingly.  Define a bigrading on $\wedge^\bullet M'$ by setting
\begin{align*}
\deg_1 x_i &= \deg_1 dx_i = 1,&  \deg_2 x_i &= \deg_2 dx_i = 0, \\
\deg_1 y_j &= -d_j,& \deg_2 y_j &= 1.
\end{align*}
Thus ``$\deg_2$'' is just ``total degree in $y$,'' which was the
grading used in section 2.  We then have 
\begin{equation}
H^{k+r}(\Omega^\bullet_{K[x,y]/K})^{(q,p)} = H^k(\wedge^\bullet
M', \bar{\partial}_h)^{(q+d_1+\cdots+d_r,p-r)}.
\end{equation}
It is then clear that
\begin{equation}
H^k(\Omega^\bullet_{K[x,y]/K})^{(q,p)} = 0 \quad\text{for all
  $k$ and $q$ if $p<r$.}
\end{equation}

The main point here is that the complex $(\wedge^\bullet
M',\bar{\partial}_h)$ is of the type studied in section 2.  Let $I$ be
the ideal of $A$ generated by the coefficients of
$df_1\wedge\cdots\wedge df_r$ relative to the basis 
\[ \{dx_{i_1}\wedge\cdots\wedge dx_{i_r} \mid 1\leq i_1<\cdots<i_r\leq
n\} \]
of~$\wedge^r M$.  We shall assume from now on that $f_1=\cdots=f_r=0$
defines a smooth complete intersection in ${\bf P}^{n-1}$, which is
equivalent to assuming that the ideal $I$ has depth~$n-r$.  Fix $q$
and apply Proposition 2.10(b) with $\omega_i$ replaced by $df_i$ and the
graded complex $((\wedge^\bullet M')^{(\bullet)},\partial)$ replaced by
$((\wedge^\bullet M')^{(q,\bullet)},\bar{\partial}_h)$ to conclude that
\[ H^k(\wedge^\bullet M',\bar{\partial}_h)^{(q,p)} = 0 \]
for $1\leq p<n-r(={\rm depth}(I))$ and $k<n-1( = {\rm
depth}(I)+r-1)$.  From (3.3) we then get the following result.

\begin{proposition}
If $0\leq k\leq n-2$ and $r+1\leq p\leq n-1$, then
\[ H^{k+r}(\Omega^\bullet_{K[x,y]/K})^{(q,p)}=0 \quad \text{for all
$q$.} \]
\end{proposition}

We need to impose a restriction on $q$ to treat the case $p\geq n$.
\begin{proposition}
If $p\geq n$ and $q\geq r-n$, then 
$H^k(\Omega^\bullet_{K[x,y]/K})^{(q,p)} = 0$ for all $k$.
\end{proposition}

{\it Proof}.  By (3.2), it suffices to prove this for $k\geq r$.  For
notational convenience we set $D=d_1+\cdots+d_r$.  By~(3.3) we are
reduced to proving 
\begin{equation}
H^k(\wedge^\bullet M',\bar{\partial}_h)^{(q,p)} =
  0\quad\text{for all $k$ when $p\geq n-r$ and $q\geq D+r-n$.}
\end{equation}
We follow the proof of Proposition 2.10(b).  For
$g_1,\ldots,g_{n-r}\in I$, put $A_i=A/(g_1,\ldots,g_i)$,
$M_i=M/(g_1,\ldots,g_i)M$, $M'_i=M'/(g_1,\ldots,g_i)M'$ for
$i=1,\ldots,n-r$.  Fix $p\geq n-r$ and $q\geq D+r-n$ and let
$\alpha_0\in(\wedge^k M')^{(q,p)}$ with
$\bar{\partial}_h(\alpha_0)=0$.  Arguing as in the proof of
Proposition 2.10(b), we construct a regular sequence
$g_1,\ldots,g_{n-r}\in I$, homogeneous elements of degrees
$e_1,\ldots,e_{n-r}$ in the grading by total degree in $x$, and elements
$\alpha_i\in(\wedge^{k-i}M'_{i-1})^{(q+e_1+\cdots+e_i,p-i)}$ such 
that
\begin{equation}
g_i\alpha_{i-1}=\bar{\partial}_h(\alpha_i)\qquad\text{in
$\wedge^{k-i+1}M'_{i-1}$}
\end{equation}
for $i=1,\ldots,n-r$.  Now consider the ring
\[ A_{n-r}=K[x]/(f_1,\ldots,f_r,g_1,\ldots,g_{n-r}) \]
and let $H_{A_{n-r}}(t)$ be its Hilbert series, i.~e.,
\[ H_{A_{n-r}}(t) = \sum_{i=0}^\infty (\dim_K A_{n-r}^{(i)})t^i, \]
where $A_{n-r}^{(i)}$ denotes the $K$-subspace of $A_{n-r}$ spanned by
polynomials homogeneous of degree $i$ in the grading by total degree
in $x$.  Since $f_1,\ldots,f_r,g_1,\ldots,g_{n-r}$ is a regular
sequence in $K[x]$, the Koszul complex on $K[x]$ defined by these
polynomials is a free resolution of $A_{n-r}$.  On exact sequences the
alternating sum of Hilbert series is zero, and since the Hilbert
series of $K[x]$ is $(1-t)^{-n}$ we get
\begin{align}
H_{A_{n-r}}(t) &= \frac{\prod_{i=1}^r(1-t^{d_i})\cdot \prod_{j=1}^{n-r}
(1-t^{e_j})}{(1-t)^n} \nonumber \\
 &= \prod_{i=1}^r (1+t+\cdots+t^{d_i-1}) \cdot \prod_{j=1}^{n-r}
(1+t+\cdots+t^{e_j-1}). 
\end{align}
Put $E=e_1+\cdots+e_{n-r}$.  We conclude that $A_{n-r}$ has no
element whose total degree in $x$ is $>D+E-n$.  But 
$\alpha_{n-r}\in(\wedge^{k-n+r} M'_{n-r-1})^{(q+E,p-n+r)}$, so from
the definition of the bigrading we see that every term
\[ x_1^{a_1}\cdots x_n^{a_n}y_1^{b_1}\cdots y_r^{b_r}\,
dx_{i_1}\wedge\cdots\wedge dx_{i_{k-n+r}} \]
appearing in $\alpha_{n-r}$ satisfies 
\begin{equation}
a_1+\cdots+a_n+(k-n+r)-\sum_{i=1}^r b_id_i= q+E. 
\end{equation}
Since $p\geq n-r\geq 1$ we have $\sum_{i=1}^r b_id_i>0$.  We also have
$q\geq D+r-n$ and $k\leq n$.  Substituting these inequalities into
(3.10) and rearranging terms gives 
\[ a_1+\cdots+a_n>D+E-n. \]
But since $A_{n-r}$ has no element whose total degree in $x$ is
$>D+E-n$, it follows that $x_1^{a_1}\cdots x_n^{a_n}=0$ in $A_{n-r}$.
Equivalently, we have in $A_{n-r-1}$ that
\[ x_1^{a_1}\cdots x_n^{a_n}\equiv 0 \pmod{g_{n-r}}. \]
Hence there exists
\[ \gamma_{n-r}\in (\wedge^{k-n+r}M'_{n-r-1})^{(q+E-e_{n-r},p-n+r)} \]
such that $\alpha_{n-r}=g_{n-r}\gamma_{n-r}$.  Substituting this in
(3.8) with $i=n-r$ and using the fact that $g_{n-r}$ is not a
zero-divisor in $A_{n-r-1}$ gives
\[ \alpha_{n-r-1}=\bar{\partial}_h(\gamma_{n-r})\qquad\text{in
$\wedge^{k-n+r+1}M'_{n-r-1}$.} \]
Suppose inductively that for some $i$, $1\leq i\leq n-r-1$, we have proved
\begin{equation}
\alpha_i=\bar{\partial}_h(\gamma_{i+1})\qquad\text{in $\wedge^{k-i} M'_i$} 
\end{equation}
for some $\gamma_{i+1}\in\wedge^{k-i-1} M'_i$.  Let
$\tilde{\gamma}_{i+1}$ be any lifting of $\gamma_{i+1}$ to
$\wedge^{k-i-1} M'_{i-1}$.  Then (3.11) implies
\[ \alpha_i=\bar{\partial}_h(\tilde{\gamma}_{i+1}) + g_i\gamma_i
\qquad\text{in $\wedge^{k-i}M'_{i-1}$} \]
for some $\gamma_i\in\wedge^{k-i} M'_{i-1}$.  Substitution into
(3.8) gives
\[ g_i\alpha_{i-1}=\bar{\partial}_h(g_i\gamma_i)\qquad\text{in
$\wedge^{k-i+1} M'_{i-1}$.} \]
Since $g_i$ is not a zero-divisor in $A_{i-1}$, we conclude that (3.11)
holds with $i$ replaced by $i-1$.  By descending induction on $i$, it
follows that (3.11) holds for $i=0$, which is the assertion of the
proposition.  

We summarize some of these observations in the following result.
\begin{proposition}
If $k\leq n+r-2$, $q\geq r-n$, and $p\neq r$, then
\[ H^k(\Omega^\bullet_{K[x,y]/K})^{(q,p)} =0. \]
\end{proposition}

{\it Proof}.  For $p\geq n$, the assertion follows from Proposition
3.6.  For $r<p<n$, the assertion follows from (3.2) when $k<r$ and
from Proposition 3.5 when $r\leq k\leq n+r-2$.  For $p<r$, the
assertion follows from (3.4). 

Equation (3.4) and Proposition 3.6 imply (1.8).  Proposition
3.12 implies (1.7) and (1.9) for $p\neq r$.  To finish the proofs of
(1.7) and (1.9), it remains only to describe the cohomology for $p=r$.
By equation (3.3) and the remark following Proposition~2.10 we have
for $k\leq n+r-2$ 
\begin{equation}
H^k(\Omega^\bullet_{K[x,y]/K})^{(q,r)} =
(\wedge^{k-2r}M)^{(q,0)}\wedge df_1\wedge\cdots\wedge df_r\wedge
dy_1\wedge\cdots\wedge dy_r. 
\end{equation}
Since $\wedge^{k-2r}M=0$ for $k<2r$, we have
\begin{equation}
H^k(\Omega^\bullet_{K[x,y]/K})^{(q,r)} = 0\quad\text{for
  all $q$ if $k<2r$.}
\end{equation}
If $r=n-1$ then $2r>n+r-2$ and there is nothing left to prove.  So
suppose $r\leq n-2$.  For $k\geq 2r$ we have
\[ (\wedge^{k-2r}M)^{(q,0)} = \bigoplus_{1\leq i_1<\cdots<i_{k-2r}\leq n}
A^{(q-k+2r)}\, dx_{i_1}\wedge\cdots\wedge dx_{i_{k-2r}}, \]
hence by (3.13)
\begin{multline}
H^k(\Omega^\bullet_{K[x,y]/K})^{(q,r)} = \\
\biggl(\sum_{1\leq
  i_1<\cdots<i_{k-2r}\leq n} A^{(q-k+2r)}\, dx_{i_1}\wedge\cdots\wedge
  dx_{i_{k-2r}}\biggr)\wedge df_1\wedge\cdots\wedge df_r\wedge
  dy_1\wedge\cdots\wedge dy_r.
\end{multline}
Now $A^{(l)}=0$ if $l<0$, so taking $q=0$ and $k>2r$ in (3.15) gives
\begin{equation}
H^k(\Omega^\bullet_{K[x,y]/K})^{(0,r)} = 0\quad\text{for
$2r<k\leq n+r-2$.}
\end{equation}
Equations (3.14) and (3.16) establish (1.7) for $p=r$.  And since
$A^{(0)}=K$, taking $q=0$ and $k=2r$ in (3.15) gives
\begin{equation}
H^{2r}(\Omega^\bullet_{K[x,y]/K},\partial)^{(0,r)} = K\cdot
  [df_1\wedge\cdots\wedge df_r\wedge dy_1\wedge\cdots\wedge dy_r],
\end{equation}
which proves (1.9) when $p=r$.  

\section{Computation of $H^{n+r-1}(\Omega^\bullet_{K[x,y]/K})$} 

Let $\theta:(\Omega_{K[x,y]/K}^k)^{(q,p)}\rightarrow
(\Omega_{K[x,y]/K}^{k-1})^{(q,p)}$ be defined by $K[x,y]$-linearity and
the formula 
\begin{multline}
\theta(dx_{i_1}\wedge\cdots\wedge dx_{i_l} \wedge
dy_{j_1}\wedge\cdots\wedge dy_{j_m})= \\
\sum_{s=1}^l (-1)^{s-1} x_{i_s} 
dx_{i_1}\wedge\cdots\wedge\widehat{dx_{i_s}}\wedge\cdots\wedge
dx_{i_l}\wedge dy_{j_1}\wedge\cdots\wedge dy_{j_m}+ \\
\sum_{t=1}^m (-1)^{l+t-1}(-d_{j_t}y_{j_t}) dx_{i_1}\wedge\cdots\wedge
dx_{i_l}\wedge dy_{j_1}\wedge\cdots\wedge\widehat{dy_{j_t}}\wedge
\cdots\wedge dy_{j_m}.
\end{multline}
One checks easily that $\theta^2=0$, $\theta(df_j)=d_jf_j$,
$\theta(dF)=0$, and
\begin{equation}
\theta(\omega_1\wedge\omega_2) = \theta(\omega_1)\wedge\omega_2 +
(-1)^m\omega_1\wedge\theta(\omega_2) 
\end{equation}
if $\omega_1$ is an $m$-form.
It follows from these latter two relations that
\begin{equation}
\theta\circ\partial + \partial\circ\theta = 0.
\end{equation}
This implies that $\theta$ induces a map
\[ \theta:H^{n+r}(\Omega^\bullet_{K[x,y]/K})^{(q,p)}\to
H^{n+r-1}(\Omega^\bullet_{K[x,y]/K})^{(q,p)}. \]
We shall prove the remaining assertions of Theorem 1.6 by studying
this induced map.  The following result is a more precise version of
assertions (1.10), (1.11), and~(1.12) of Theorem 1.6.

\begin{theorem}
Suppose that $r<n$ and that the equations $f_1=\dots=f_r=0$ define a
smooth complete intersection $X$ in ${\bf P}^{n-1}$.  \\
{\bf (a)} Assume $r<n-1$.  If $d_1\cdots d_r\neq 0$ in $K$, then 
\[ \theta:H^{n+r}(\Omega^\bullet_{K[x,y]/K})^{(0,p)}\to
H^{n+r-1}(\Omega^\bullet_{K[x,y]/K})^{(0,p)} \]
is an isomorphism for
all $p$.  If $d_1\cdots d_r=0$ in $K$, it is an isomorphism  
for all $p$ except in the following three cases: if $n+r$ is odd and
$p=(n+r-1)/2$, it has a one-dimensional kernel and cokernel; if
$n+r$ is even and $p=(n+r)/2$, it is surjective and has a
one-dimensional kernel; and if $n+r$ is even and $p=(n+r)/2 - 1$, 
it is injective and has a one-dimensional cokernel. \\
{\bf (b)} Assume $r=n-1$.  Then
\[ \theta:H^{n+r}(\Omega^\bullet_{K[x,y]/K})^{(0,p)}\to 
H^{n+r-1}(\Omega^\bullet_{K[x,y]/K})^{(0,p)} \]
is an isomorphism for all $p$ except $p=r$, in which case it is
injective and has a one-dimensional cokernel.
\end{theorem}

Note that the complex 
\begin{equation}
0\rightarrow \Omega^{n+r}_{K[x,y]/K}\xrightarrow{\theta}\cdots
\xrightarrow{\theta}\Omega^0_{K[x,y]/K} \rightarrow 0 
\end{equation}
is isomorphic to the Koszul complex on $K[x,y]$ defined by the elements
\[ x_1,\ldots,x_n,-d_1y_1,\ldots,-d_ry_r. \]
When $d_1\cdots d_r\neq 0$ in $K$, these elements form a regular sequence 
so this complex is exact except at the right-hand term and the
following result is clear (with no restriction on $q$).  It is
somewhat surprising that it holds without any restriction on the
characteristic of $K$. 

\begin{proposition}
For all $q\geq 0$, the sequence
\[ 0\rightarrow (\Omega^{n+r}_{K[x,y]/K})^{(q,p)}\xrightarrow{\theta}\cdots
\xrightarrow{\theta}(\Omega^0_{K[x,y]/K})^{(q,p)}\to \left.
\begin{cases} K & \text{if $(q,p) = (0,0)$} \\ 0 & \text{otherwise}
\end{cases} \right\} \rightarrow 0 \]
is exact.
\end{proposition}

{\it Proof}.  Suppose $d_1\cdots d_s\neq 0$ in $K$ but $d_j=0$ in $K$
for $j=s+1,\dots,r$.  The complex (4.5) is then 
isomorphic to the Koszul complex on $K[x,y]$ defined by the elements
\[ x_1,\dots,x_n,-d_1y_1,\dots,-d_sy_s,0,\dots,0\;\text{($r-s$
times).} \]
The elements $x_1,\dots,x_n,-d_1y_1,\dots,-d_sy_s$ form a regular
sequence on $K[x,y]$.  It is then straightforward to calculate that
the quotient
\[ \frac{\ker(\theta:\Omega^k_{K[x,y]/K}\to\Omega^{k-1}_{K[x,y]/K})}
{{\rm im}(\theta:\Omega^{k+1}_{K[x,y]/K}\to\Omega^k_{K[x,y]/K})}, \]
the $k$-th homology of the complex (4.5), vanishes for $k>r-s$ and
for $k\leq r-s$  is isomorphic to 
\begin{equation}
\bigoplus_{s+1\leq j_1<\dots<j_k\leq r}K[y_{s+1},\dots,y_r]\,
dy_{j_1}\wedge\cdots\wedge dy_{j_k}
\end{equation}
with the induced bigrading.  But since $\deg_1 y_j$ and $\deg_1 dy_j$
are negative, we have for $q\geq 0$ that
\[ \biggl(\bigoplus_{s+1\leq j_1<\dots<j_k\leq r}K[y_{s+1},\dots,y_r]\,
dy_{j_1}\wedge\cdots\wedge dy_{j_k}\biggr)^{(q,p)} = \begin{cases} K &
\text{if $k=q=p=0$}, \\ 0 & \text{otherwise.} \end{cases} \]
This establishes the proposition.

Proposition 4.6 allows us to construct a short exact sequence of
complexes involving
$(\Omega^\bullet_{K[x,y]/K},\partial)^{(0,\bullet)}$.  For $i\geq 0$, put
\begin{equation}
\widetilde{\Omega}^i_{K[x,y]/K}=
\theta(\Omega^{i+1}_{K[x,y]/K})\subseteq\Omega^i_{K[x,y]/K}.  
\end{equation}
Equation (4.3) implies that
$\partial(\widetilde{\Omega}^i_{K[x,y]/K})\subseteq 
\widetilde{\Omega}^{i+1}_{K[x,y]/K}$,
thus $(\widetilde{\Omega}^\bullet_{K[x,y]/K},\partial)$ is a
subcomplex of $(\Omega^\bullet_{K[x,y]/K},\partial)$.    

We define a related complex $\widehat{\Omega}^\bullet_{K[x,y]/K}$ as
follows.  Let 
\begin{equation}
\widehat{\Omega}^0_{K[x,y]/K}=
\Omega^0_{K[x,y]/K}/\widetilde{\Omega}^0_{K[x,y]/K} 
\end{equation}
and let $\widehat{\Omega}^i_{K[x,y]/K} =
\widetilde{\Omega}^{i-1}_{K[x,y]/K}$ 
for $i\geq 1$.  We define the boundary map
$\widehat{\Omega}^i_{K[x,y]/K} \rightarrow
\widehat{\Omega}^{i+1}_{K[x,y]/K}$ to be zero if $i=0$ and $-\partial$
if $i\geq 1$.  Thus 
\begin{equation}
H^0(\widehat{\Omega}^\bullet_{K[x,y]/K}) =
\Omega^0_{K[x,y]/K}/\widetilde{\Omega}^0_{K[x,y]/K}
\end{equation}
and $H^i(\widehat{\Omega}^\bullet_{K[x,y]/K}) =
H^{i-1}(\widetilde{\Omega}^\bullet_{K[x,y]/K})$ for $i\geq 1$.  Define maps
$\Omega^i_{K[x,y]/K} \rightarrow \widehat{\Omega}^i_{K[x,y]/K}$ as
follows.  For $i=0$, take the map that induces the isomorphism (4.9) and for
$i\geq 1$ take the map $\theta:\Omega^i_{K[x,y]/K} \rightarrow
\widehat{\Omega}^i_{K[x,y]/K}$.  It follows from Proposition~4.6 that these
maps define short exact sequences of complexes for all $q\geq 0$:
\begin{equation}
0\rightarrow (\widetilde{\Omega}^\bullet_{K[x,y]/K})^{(q,\bullet)}\rightarrow
(\Omega^\bullet_{K[x,y]/K})^{(q,\bullet)} \xrightarrow{\theta}
(\widehat{\Omega}^\bullet_{K[x,y]/K})^{(q,\bullet)} \rightarrow 0.
\end{equation}
Note that these maps respect the bigrading defined earlier.
We thus get exact sequences of cohomology groups
\begin{multline}
\cdots\rightarrow H^i(\widetilde{\Omega}^\bullet_{K[x,y]/K})^{(q,p)}
\rightarrow H^i(\Omega^\bullet_{K[x,y]/K})^{(q,p)} \xrightarrow{\theta} \\
H^{i-1}(\widetilde{\Omega}^\bullet_{K[x,y]/K})^{(q,p)}\xrightarrow{\delta}
H^{i+1}(\widetilde{\Omega}^\bullet_{K[x,y]/K})^{(q,p+1)} \rightarrow \cdots.
\end{multline}
Note that the connecting homomorphism $\delta$ increases $\deg_2$ by 1, i.~e., 
\[ \delta(H^{i-1}(\widetilde{\Omega}^\bullet_{K[x,y]/K})^{(q,p)}) \subseteq
H^{i+1}(\widetilde{\Omega}^\bullet_{K[x,y]/K})^{(q,p+1)}. \]

The exact sequence (4.12) shows that $\theta$ induces an isomorphism
\begin{equation}
H^{n+r}(\Omega^\bullet_{K[x,y]/K})^{(0,p)} \cong
H^{n+r-1}(\widetilde{\Omega}^\bullet_{K[x,y]/K})^{(0,p)}
\quad\text{for all $p$.}
\end{equation}
The map $\theta:H^{n+r}(\Omega^\bullet_{K[x,y]/K})^{(0,p)}\to
H^{n+r-1}(\Omega^\bullet_{K[x,y]/K})^{(0,p)}$ of Theorem 4.4 factors
through this isomorphism as
\[ H^{n+r}(\Omega^\bullet_{K[x,y]/K})^{(0,p)}\xrightarrow{\theta}
H^{n+r-1}(\widetilde{\Omega}^\bullet_{K[x,y]/K})^{(0,p)}\to
H^{n+r-1}(\Omega^\bullet_{K[x,y]/K})^{(0,p)}, \]
where the second map is induced by the inclusion
\[ (\widetilde{\Omega}^\bullet_{K[x,y]/K},\partial)^{(0,p)}
\hookrightarrow (\Omega^\bullet_{K[x,y]/K},\partial)^{(0,p)}. \]
Thus to prove Theorem 4.4, it suffices to prove the asserted
properties for the map
\begin{equation}
H^{n+r-1}(\widetilde{\Omega}^\bullet_{K[x,y]/K})^{(0,p)}\to 
H^{n+r-1}(\Omega^\bullet_{K[x,y]/K})^{(0,p)}.
\end{equation}
We shall accomplish this by computing the cohomology of all the terms
of (4.12) when $q=0$.

By (1.7), if $i<2r$ then $H^i(\Omega^\bullet_{K[x,y]/K})^{(0,p)} = 0$
for all $p$.  Using this fact in (4.12) shows that the connecting
homomorphism $\delta$ gives isomorphisms
\begin{equation}
H^0(\widehat{\Omega}^\bullet_{K[x,y]/K})^{(0,p)}\cong
H^1(\widetilde{\Omega}^\bullet_{K[x,y]/K})^{(0,p+1)} 
\end{equation}
and
\begin{equation}
H^i(\widetilde{\Omega}^\bullet_{K[x,y]/K})^{(0,p)} \cong
H^{i+2}(\widetilde{\Omega}^\bullet_{K[x,y]/K})^{(0,p+1)}  
\end{equation}
for $i=0,1,\dots,2r-3$ and all $p$.  From (4.10) and Proposition 4.6 we have
\begin{equation}
H^0(\widehat{\Omega}^\bullet_{K[x,y]/K})^{(0,p)} = \begin{cases} K &
\text{if $p=0$,} \\ 0 & \text{otherwise} \end{cases}
\end{equation}
and from (1.7) and (4.12) we have
\begin{equation}
H^0(\widetilde{\Omega}^\bullet_{K[x,y]/K})^{(0,p)}=0 \quad\text{for
all $p$.}
\end{equation}
Using (4.15), (4.17), and (4.18), it now follows inductively from (4.16) that
\begin{equation}
H^{2k}(\widetilde{\Omega}^\bullet_{K[x,y]/K})^{(0,p)} = 0 \quad \text{for
$0\leq k<r$ and all $p$} 
\end{equation}
and 
\begin{equation}
H^{2k-1}(\widetilde{\Omega}^\bullet_{K[x,y]/K})^{(0,p)} \cong \left.
\begin{cases} K & \text{if $p=k$} \\ 0 & \text{otherwise} \end{cases}
\right\} \quad\text{for $1\leq k\leq r$.}
\end{equation}

It is useful to specify a basis $[\eta_k]$ for
$H^{2k-1}(\widetilde{\Omega}^\bullet_{K[x,y]/K})^{(0,k)}$.  By (4.17),
the class $[1]$ is a basis for
$H^0(\widehat{\Omega}^\bullet_{K[x,y]/K})^{(0,0)}$ so we define
$\eta_0=1$.  The isomorphism~(4.15) given by $\delta$ sends $[1]$ to
$[dF]$, so define $\eta_1 = dF$.  Now let $2\leq k\leq r-1$ and
suppose that $\eta_{k-1}\in
(\widetilde{\Omega}^{2k-3}_{K[x,y]/K})^{(0,k-1)}$ has been chosen such that
$[\eta_{k-1}]$ is a basis for
$H^{2k-3}(\widetilde{\Omega}^\bullet_{K[x,y]/K})^{(0,k-1)}$.  Choose
$\zeta_{k-1}\in(\Omega^{2k-2}_{K[x,y]/K})^{(0,k-1)}$ such that
\begin{equation}
\theta(\zeta_{k-1}) = \eta_{k-1}
\end{equation}
and define
\begin{equation}
\eta_k = dF\wedge\zeta_{k-1}\in
(\widetilde{\Omega}^{2k-1}_{K[x,y]/K})^{(0,k)}.
\end{equation}
The definition of $\delta$ shows that $[\eta_k]$ is the image of
$[\eta_{k-1}]$ under the isomorphism (4.16), hence $[\eta_k]$ is a
basis for $H^{2k-1}(\widetilde{\Omega}^\bullet_{K[x,y]/K})^{(0,k)}$.

The following result is the key to calculating the
$H^i(\widetilde{\Omega}^\bullet_{K[x,y]/K})$ for $i\geq 2r$.  For
$k=1,\dots,r$, let $\xi_k\in (\Omega^{2k}_{K[x,y]/K})^{(0,k)}$ be
defined by
\begin{equation}
\xi_k=\sum_{1\leq i_1<\cdots<i_k\leq r}
\biggl(\prod_{i\not\in\{i_1,\dots,i_k\}} d_i\biggr)\,
df_{i_1}\wedge\cdots\wedge df_{i_k}\wedge dy_{i_1}\wedge\cdots\wedge
dy_{i_k}.
\end{equation}

\begin{proposition}
Let $r<n-1$.  Relative to the bases $[\xi_r]$ for
$H^{2r}(\Omega^\bullet_{K[x,y]/K})^{(0,r)}$ and $[\eta_r]$ for
$H^{2r-1}(\widetilde{\Omega}^\bullet_{K[x,y]/K})^{(0,r)}$, the map
\[ \theta:H^{2r}(\Omega^\bullet_{K[x,y]/K})^{(0,r)}\to
H^{2r-1}(\widetilde{\Omega}^\bullet_{K[x,y]/K})^{(0,r)} \]
is multiplication by $(-1)^{r(r-1)/2} d_1\cdots d_r$.
\end{proposition}

{\it Proof}.  We prove inductively that for $k=1,\dots,r$,
\begin{equation}
\theta(\xi_k) = (-1)^{k(k-1)/2}(d_1\cdots d_r)\eta_k +
dF\wedge\theta(\tau_k) 
\end{equation}
for some $\tau_k\in (\Omega^{2k-1}_{K[x,y]/K})^{(0,k-1)}$.  The
assertion of the proposition follows by taking $k=r$ in (4.25).  For
$k=1$, a straightforward calculation shows that
\[ \theta(\xi_1) = (d_1\cdots d_r)dF = (d_1\cdots d_r)\eta_1, \]
so suppose (4.25) holds for some $k$, $1\leq k<r$.  A straightforward
calculation shows that
\begin{equation}
\theta(\xi_{k+1}) = (-1)^k dF\wedge\xi_k.
\end{equation}
As in (4.21), choose $\zeta_k$ so that $\theta(\zeta_k) = \eta_k$.
Substitution into (4.25) then gives
\[ \theta(\xi_k) = \theta((-1)^{k(k-1)/2}(d_1\cdots d_r)\zeta_k +
dF\wedge\tau_k) \] 
(since $\theta(dF) = 0$), hence by Proposition~4.6 there exists
$\tau_{k+1}$ such that
\[ \xi_k = (-1)^{k(k-1)/2}(d_1\cdots d_r)\zeta_k + dF\wedge\tau_k +
\theta((-1)^k\tau_{k+1}). \]
Substitution into (4.26) now gives
\[ \theta(\xi_{k+1}) = (-1)^{k(k+1)/2}(d_1\cdots d_r)dF\wedge\zeta_k +
dF\wedge\theta(\tau_{k+1}). \]
Since $dF\wedge\zeta_k = \eta_{k+1}$ by (4.22), this is just (4.25)
with $k$ replaced by $k+1$.

\begin{corollary}
Suppose $r<n-1$ and $d_1\cdots d_r\neq 0$ in $K$.  Then the map
\[ \theta:H^{2r}(\Omega^\bullet_{K[x,y]/K})^{(0,p)}\to 
H^{2r-1}(\widetilde{\Omega}^\bullet_{K[x,y]/K})^{(0,p)} \] 
is an isomorphism for all $p$.
\end{corollary}

{\it Proof}.  By (1.9) and (4.20), both cohomology groups vanish if
$p\neq r$.  If $p=r$, $\theta$~is an isomorphism by
Proposition~4.24. 

\begin{lemma} Suppose $d_1\cdots d_r\neq 0$ in $K$.  Then
\[ H^i(\widetilde{\Omega}^\bullet_{K[x,y]/K})^{(0,p)} =
0\quad\text{for $2r\leq i\leq n+r-2$ and all $p$.} \]
\end{lemma}

{\it Proof}.  For $r=n-1$ there is nothing to prove (since $2r>n+r-2$
in that case), so assume $r<n-1$.  Using (4.19), (4.12) gives an
exact sequence
\[ 0\to H^{2r}(\widetilde{\Omega}^\bullet_{K[x,y]/K})^{(0,p)}\to 
H^{2r}(\Omega^\bullet_{K[x,y]/K})^{(0,p)} \xrightarrow{\theta}
H^{2r-1}(\widetilde{\Omega}^\bullet_{K[x,y]/K})^{(0,p)}. \]
It then follows from Corollary 4.27 that
\[ H^{2r}(\widetilde{\Omega}^\bullet_{K[x,y]/K})^{(0,p)} =
0\quad\text{for all $p$.} \]
If $r=n-2$, then $2r=n+r-2$ and we are done.  So assume also $r<
n-2$.  Then $2r+1<n+r-1$, so
$H^{2r+1}(\Omega^\bullet_{K[x,y]/K})^{(0,p+1)} = 0$ by (1.7), and (4.12)
gives an exact sequence 
\[ H^{2r}(\Omega^\bullet_{K[x,y]/K})^{(0,p)} \xrightarrow{\theta}
H^{2r-1}(\widetilde{\Omega}^\bullet_{K[x,y]/K})^{(0,p)}
\xrightarrow{\delta}
H^{2r+1}(\widetilde{\Omega}^\bullet_{K[x,y]/K})^{(0,p+1)} \to 0. \]
It now follows from Corollary 4.27 that
\[ H^{2r+1}(\widetilde{\Omega}^\bullet_{K[x,y]/K})^{(0,p)} =
0\quad\text{for all $p$.} \]

Assume now that for some $i$, $2r<i<n+r-2$, we have proved
\begin{equation}
H^{i-1}(\widetilde{\Omega}^\bullet_{K[x,y]/K})^{(0,p)} =
H^i(\widetilde{\Omega}^\bullet_{K[x,y]/K})^{(0,p)}=0\quad\text{for all
  $p$.} 
\end{equation}
Using (1.7) in the exact sequence (4.12) gives
\begin{align*}
H^{i-1}(\widetilde{\Omega}^\bullet_{K[x,y]/K})^{(0,p)} &\cong
H^{i+1}(\widetilde{\Omega}^\bullet_{K[x,y]/K})^{(0,p+1)} \quad\text{for
all $p$} \\
 &=0 \quad\text{by (4.29)}.
\end{align*}
The assertion of the lemma now follows by induction on $i$.

We can now prove Theorem 4.4 in the case where $d_1\cdots d_r\neq 0$
in $K$.  Suppose first $r<n-2$.  Using Lemma~4.28 with $i=n+r-3,n+r-2$
in (4.12) shows that the map (4.14) is an isomorphism for all $p$.  If
$r=n-2$, using Lemma 4.28 with $i=n+r-2$ in (4.12) gives an exact sequence
\begin{multline*}
H^{2r-1}(\widetilde{\Omega}^\bullet_{K[x,y]/K})^{(0,p-1)}
\xrightarrow{\delta}
H^{n+r-1}(\widetilde{\Omega}^\bullet_{K[x,y]/K})^{(0,p)} \to \\
H^{n+r-1}(\Omega^\bullet_{K[x,y]/K})^{(0,p)} \to 0.
\end{multline*}
If $p\neq r+1$, then $\delta$ is the zero map by (4.20).  If $p=r+1$,
then by Proposition~4.24 the image of $\delta$ is spanned by
\[ \delta([\theta(\xi_r)]) = [\partial(\xi_r)] = 0, \]
so $\delta$ is the zero map in this case also.  Thus (4.14) is an
isomorphism for $r=n-2$ also.   If $r=n-1$, using (4.19) in (4.12) 
gives an exact sequence 
\begin{multline}
0\to H^{n+r-1}(\widetilde{\Omega}^\bullet_{K[x,y]/K})^{(0,p)}\to
H^{n+r-1}(\Omega^\bullet_{K[x,y]/K})^{(0,p)} \xrightarrow{\theta} \\
H^{2r-1}(\widetilde{\Omega}^\bullet_{K[x,y]/K})^{(0,p)}\to 0.
\end{multline}
Then by (4.20), the map (4.14) is an isomorphism for $p\neq r$ and is
injective and has a one-dimensional cokernel for $p=r$.

When $d_1\cdots d_r = 0$ in $K$, Proposition 4.24 gives the following.
\begin{corollary}
Suppose $r<n-1$ and $d_1\cdots d_r=0$ in $K$.  Then the map
\[ \theta:H^{2r}(\Omega^\bullet_{K[x,y]/K})^{(0,p)}\to
H^{2r-1}(\widetilde{\Omega}^{\bullet}_{K[x,y]/K})^{(0,p)} \]
is the zero map for all $p$.
\end{corollary}

This leads to the following result.
\begin{lemma}
Suppose $d_1\cdots d_r=0$ in $K$ and $2r\leq i\leq n+r-2$.
If $i$ is even, then
\[ \dim_K H^i(\widetilde{\Omega}^\bullet_{K[x,y]/K})^{(0,p)} =
\begin{cases} 1 & \text{if $p=i/2$} \\ 0 & \text{otherwise}
\end{cases} \]
and if $i$ is odd, then
\[ \dim_K H^i(\widetilde{\Omega}^\bullet_{K[x,y]/K})^{(0,p)} =
\begin{cases} 1 & \text{if $p=(i+1)/2$} \\ 0 & \text{otherwise.}
\end{cases} \]
\end{lemma}

{\it Proof}.  For $r=n-1$ there is nothing to prove, so suppose $r\leq
n-2$.  Using (4.19) and Corollary 4.31 in (4.12) gives isomorphisms for all $p$
\[ H^{2r}(\widetilde{\Omega}_{K[x,y]/K}^\bullet)^{(0,p)} \cong 
H^{2r}(\Omega_{K[x,y]/K}^\bullet)^{(0,p)}. \]
The assertion of the lemma now follows for $i=2r$ by (1.9).  If
$r=n-2$ there is nothing left to prove, so suppose $r<n-2$.  Using
(1.7) and Corollary~4.31 in (4.12) gives isomorphisms for all $p$
\[ H^{2r-1}(\widetilde{\Omega}_{K[x,y]/K}^\bullet)^{(0,p)} \cong 
H^{2r+1}(\widetilde{\Omega}_{K[x,y]/K}^\bullet)^{(0,p+1)}. \] 
The assertion of the lemma now follows for $i=2r+1$ by (4.20).  If
$r=n-3$, there is nothing left to prove so suppose $r<n-3$.  
Suppose inductively the lemma is true for some $i$, $2r\leq 
i\leq n+r-4$.  By (1.7) we have
\[ H^{i+1}(\Omega^\bullet_{K[x,y]/K}) =
H^{i+2}(\Omega^\bullet_{K[x,y]/K}) = 0, \]
so (4.12) gives isomorphisms for all $p$
\[ H^i(\widetilde{\Omega}_{K[x,y]/K}^\bullet)^{(0,p)} \cong 
H^{i+2}(\widetilde{\Omega}_{K[x,y]/K}^\bullet)^{(0,p+1)}. \] 
The assertion of the lemma now follows for $i+2$, and by induction on
$i$ the proof is complete.

We can now prove Theorem 4.4 when $d_1\cdots d_r=0$ in $K$.
If $r<n-1$, then the map
$\theta:H^{n+r-2}(\Omega^\bullet_{K[x,y]/K})^{(0,p)}\to  
H^{n+r-3}(\widetilde{\Omega}^\bullet_{K[x,y]/K})^{(0,p)}$ is the zero
map for all $p$ (use (1.7) if $r<n-2$ and use Corollary 4.31 if
$r=n-2$), so (4.12) gives an exact sequence
\begin{multline}
0\to H^{n+r-3}(\widetilde{\Omega}^\bullet_{K[x,y]/K})^{(0,p-1)}
\xrightarrow{\delta}
H^{n+r-1}(\widetilde{\Omega}^\bullet_{K[x,y]/K})^{(0,p)} \to \\
H^{n+r-1}(\Omega^\bullet_{K[x,y]/K})^{(0,p)}\xrightarrow{\theta}
H^{n+r-2}(\widetilde{\Omega}^\bullet_{K[x,y]/K})^{(0,p)}\to 0.
\end{multline}
Applying Lemma 4.32 to this exact sequence shows that (4.14) is an
isomorphism for all $p$ except in three cases.  If $n+r$ is odd and
$p=(n+r-1)/2$, the map (4.14) has a one-dimensional kernel and cokernel. 
If $n+r$ is even and $p=(n+r)/2$, the map (4.14) is surjective and has
a one-dimensional kernel, while if $p=(n+r)/2-1$, the map (4.14) is
injective and has a one-dimensional cokernel.
If $r=n-1$ (so that $2r=n+r-1$), then using (4.19) in (4.12) gives
short exact sequences for all~$p$
\begin{multline}
0\to H^{n+r-1}(\widetilde{\Omega}^\bullet_{K[x,y]/K})^{(0,p)} \to
H^{n+r-1}(\Omega^\bullet_{K[x,y]/K})^{(0,p)} \xrightarrow{\theta} \\
H^{n+r-2}(\widetilde{\Omega}^\bullet_{K[x,y]/K})^{(0,p)} \to 0.
\end{multline}
Now $n+r-2=2r-1$, so we have by (4.20) that
\[ \dim_K H^{n+r-2}(\widetilde{\Omega}^\bullet_{K[x,y]/K})^{(0,p)} =
\begin{cases} 1 & \text{if $p=r$} \\ 0 & \text{otherwise.}
\end{cases} \]
It then follows from (4.34) that (4.14) is an isomorphism for all $p$
except $p=r$, in which case it is injective and has a one-dimensional
cokernel.  This completes the proof of Theorem 4.4 (and hence the
proof of Theorem 1.6).  

\section{Hilbert series of $H^{n+r}(\Omega^\bullet_{K[x,y]/K})^{(0,\bullet)}$}

In this section we compute the Hilbert series of
$H^{n+r}(\Omega^\bullet_{K[x,y]/K})^{(0,\bullet)}$, i.~e., the series
\[ \sum_{p=0}^\infty (\dim_K
H^{n+r}(\Omega^\bullet_{K[x,y]/K})^{(0,p)})t^p, \]
which, by (1.8), is a polynomial of degree $\leq n-1$ divisible by $t^r$.

A basis for $(K[x,y]\,dx_{i_1}\cdots dx_{i_l}dy_{j_1}\cdots
dy_{j_m})^{(0,p)}$ is given by the forms
\begin{equation}
x_1^{a_1}\cdots x_n^{a_n}y_1^{b_1}\cdots y_r^{b_r}\,dx_{i_1}\cdots
dx_{i_l}dy_{j_1}\cdots dy_{j_m}
\end{equation}
with
\begin{equation}
a_1+\cdots+a_n = b_1d_1+\cdots+b_rd_r+d_{j_1}+\cdots+d_{j_m}-l
\end{equation}
and
\begin{equation}
b_1+\cdots+b_r+m=p.
\end{equation}

Define polynomials $p_l(b_1,\dots,b_r)\in{\bf Q}[b_1,\dots,b_r]$ by
\begin{equation}
p_l(b_1,\dots,b_r) = \frac{1}{(n-1)!} \prod_{j=1}^{n-1}
(b_1d_1+\cdots+b_rd_r-l+j). 
\end{equation}
For fixed $b_1,\dots,b_r,j_1,\dots,j_m,l$, the number of sequences
$a_1,\dots,a_n$ of nonnegative integers satisfying (5.2) is given by
the binomial coefficient
\[ \binom{b_1d_1+\cdots+b_rd_r+d_{j_1}+\cdots+d_{j_m}-l+n-1}{n-1}, \]
which is understood to be 0 when
\[ b_1d_1+\cdots+b_rd_r+d_{j_1}+\cdots+d_{j_m}-l < 0. \]
In terms of the polynomial (5.4), this equals
\begin{multline}
p_l(b_1,\dots,b_{j_1}+1,\dots,b_{j_m}+1,\dots,b_r) + \\
\begin{cases} (-1)^n  & \text{if $b_i=0$ for all $i$, $m=0$, $l=n$}, \\
0 & \text{otherwise}.
\end{cases}
\end{multline}

Consider the series
\begin{multline}
H_l(j_1,\dots,j_m;t_1,\dots,t_r) = \\
\sum_{b_1,\dots,b_r=0}^{\infty}
p_l(b_1,\dots,b_{j_1}+1,\dots,b_{j_m}+1,\dots,b_r)t_1^{b_1}\cdots
t_r^{b_r} t_{j_1}\cdots t_{j_m}.
\end{multline}
It follows from (5.5) that
\begin{multline}
\sum_{p=0}^{\infty} (\dim_K (K[x,y]\,dx_{i_1}\cdots dx_{i_l}
dy_{j_1}\cdots dy_{j_m})^{(0,p)}) t^p = \\
H_l(j_1,\dots,j_m;t,\dots,t) + \begin{cases} 
(-1)^n & \text{if $l=n$, $m=0$}, \\ 0 & \text{otherwise},
\end{cases}
\end{multline}
hence the Hilbert series of the complex
$(\Omega^\bullet_{K[x,y]/K},\partial)^{(0,\bullet)}$ is
\begin{multline}
(-1)^{n+r}t^r + \\ \sum_{l=0}^n \sum_{m=0}^r \sum_{1\leq
j_1<\cdots<j_m\leq r} (-1)^{n+r-l-m}\binom{n}{l}t^{n+r-l-m}
H_l(j_1,\dots,j_m;t,\dots,t),
\end{multline}
i.~e., the coefficient of $t^p$ in this series is the alternating sum
of the dimensions of the terms in the sequence
\[ 0 \to (\Omega^0_{K[x,y]/K})^{(0,p-n-r)}\xrightarrow{\partial}
\cdots \xrightarrow{\partial} (\Omega^{n+r}_{K[x,y]/K})^{(0,p)} \to
0. \]

To simplify (5.8) we begin by observing that $p_l(b_1,\dots,b_r)$ is a
polynomial of degree $n-1$, say,
\begin{equation}
p_l(b_1,\dots,b_r) = \sum_{e_1+\cdots+e_r\leq n-1} a^{(l)}_{e_1\dots
e_r} b_1^{e_1}\cdots b_r^{e_r}.
\end{equation}
The coefficients $a^{(l)}_{e_1\dots e_r}$ can be computed explicitly
from (5.4) (for simplicity we set $E=e_1+\cdots+e_r$): 
\begin{equation}
a^{(l)}_{e_1\dots e_r} = 
\frac{(-1)^{n-1-E}E!}{(n-1)!\,e_1!\cdots e_r!} s_{n-1-E}(l-(n-1),\dots,
l-1) d_1^{e_1}\cdots d_r^{e_r},
\end{equation}
where $s_i$ denotes the $i$-th elementary symmetric function in $n-1$
variables.  From (5.6) we get
\begin{multline}
H_l(j_1,\dots,j_m;t_1,\dots,t_r) = \\
 \sum_{E\leq n-1}a^{(l)}_{e_1\dots e_r}
\sum_{b_1,\dots,b_r=0}^\infty b_1^{e_1}\cdots
(b_{j_1}+1)^{e_{j_1}}\cdots (b_{j_m}+1)^{e_{j_m}}\cdots b_r^{e_r}
t_1^{b_1}\cdots t_r^{b_r}t_{j_1}\cdots t_{j_m}.
\end{multline}
Note that 
\begin{align*}
\sum_{b=0}^\infty b^et^b &= \bigg(t\frac{d}{dt}\bigg)^e\frac{1}{1-t} \\
\sum_{b=0}^\infty (b+1)^et^{b+1} &= \bigg(t\frac{d}{dt}\bigg)^e\frac{t}{1-t}.
\end{align*}
Define polynomials $p_e(t)$, $\tilde{p}_e(t)$, by
\begin{align}
\frac{p_e(t)}{(1-t)^{e+1}} &= \bigg(t\frac{d}{dt}\bigg)^e\frac{1}{1-t} \\
\frac{\tilde{p}_e(t)}{(1-t)^{e+1}} &= \bigg(t\frac{d}{dt}\bigg)^e\frac{t}{1-t}.
\end{align}
From (5.11) we then get
\begin{multline}
H_l(j_1,\dots,j_m;t,\dots,t) =  \\
\sum_{E\leq n-1}
a^{(l)}_{e_1\dots e_r}\frac{p_{e_1}(t)\cdots \tilde{p}_{e_{j_1}}(t)
\cdots \tilde{p}_{e_{j_m}}(t)\cdots p_{e_r}(t)}{(1-t)^{E+r}}.
\end{multline}

Since $1/(1-t) = 1+(t/(1-t))$, we have
\begin{equation}
p_e(t) = \tilde{p}_e(t) \quad\text{if $e>0$},
\end{equation}
while
\begin{equation}
p_0(t) = 1 \quad\text{and}\quad \tilde{p}_0(t) = t.
\end{equation}
For fixed $e_1,\cdots,e_r$, we claim that
\begin{multline}
\sum_{m=0}^r \sum_{1\leq j_1<\cdots<j_m\leq r} (-1)^{r-m}t^{r-m}
p_{e_1}(t)\cdots \tilde{p}_{e_{j_1}}(t) \cdots
\tilde{p}_{e_{j_m}}(t)\cdots p_{e_r}(t) = \\
\begin{cases}
(1-t)^r p_{e_1}(t)\cdots p_{e_r}(t) & \text{if $e_i\geq 1$ for all
$i$}, \\
0 & \text{if $e_i=0$ for some $i$}.
\end{cases}
\end{multline}
In the first case, it follows from (5.15) that the left-hand side of
(5.17) equals
\[ p_{e_1}(t)\cdots
p_{e_r}(t)\sum_{m=0}^r\binom{r}{m}(-1)^{r-m}t^{r-m}, \]
which clearly equals the right-hand side of (5.17) in that case.
In the second case, suppose, to fix ideas, that $e_1=0$.  We use 
(5.16) to break the inner sum in (5.17) into two parts, the first a
sum of those terms where $j_1>1$, the second a sum of those terms
where $j_1=1$:  
\begin{multline*}
\sum_{m=0}^r \bigg(\sum_{2\leq j_1<\dots<j_m\leq r} (-1)^{r-m}t^{r-m}
p_{e_2}(t)\cdots \tilde{p}_{e_{j_1}}(t) \cdots 
\tilde{p}_{e_{j_m}}(t)\cdots p_{e_r}(t) + \\
\sum_{2\leq j_2<\dots<j_m\leq r}(-1)^{r-m}t^{r-m+1} p_{e_2}(t)\cdots
\tilde{p}_{e_{j_2}}(t) \cdots \tilde{p}_{e_{j_m}}(t)\cdots p_{e_r}(t)
\bigg). 
\end{multline*}
This may be rewritten as
\begin{multline*}
\sum_{m=0}^{r-1} \sum_{2\leq j_1<\dots<j_m\leq r} (-1)^{r-m}t^{r-m}
p_{e_2}(t)\cdots \tilde{p}_{e_{j_1}}(t) \cdots 
\tilde{p}_{e_{j_m}}(t)\cdots p_{e_r}(t) + \\
\sum_{m=1}^r \sum_{2\leq j_2<\dots<j_m\leq r}(-1)^{r-m}t^{r-m+1}
p_{e_2}(t)\cdots \tilde{p}_{e_{j_2}}(t) \cdots
\tilde{p}_{e_{j_m}}(t)\cdots p_{e_r}(t).
\end{multline*}
Shifting the index $m$ down by 1 in the second double sum, one sees
that the second double sum is the negative of the first, which proves
(5.17) in the second case.

Substituting (5.14) in (5.8) and using (5.17), it follows that the
Hilbert series of $(\Omega^\bullet_{K[x,y]/K},\partial)^{(0,\bullet)}$
is
\begin{equation}
(-1)^{n+r}t^r + \sum_{l=0}^n (-1)^{n-l}\binom{n}{l}t^{n-l}
\sum_{\substack{E\leq n-1\\ e_i\geq 1\text{ for all $i$}}}
a^{(l)}_{e_1\dots e_r} \frac{p_{e_1}(t)\cdots
p_{e_r}(t)}{(1-t)^E}.
\end{equation}

Let $H(t)$ be as defined in (1.13).  Using Theorem 1.6, one can express
the Hilbert series of
$H^{n+r}(\Omega^\bullet_{K[x,y]/K})^{(0,\bullet)}$ and
$H^{n+r-1}(\Omega^\bullet_{K[x,y]/K})^{(0,\bullet)}$ in terms of
$H(t)$.  One then calculates that in all cases, the Hilbert
series of $(\Omega^\bullet_{K[x,y]/K},\partial)^{(0,\bullet)}$ equals
\begin{equation}
(1-t)H(t) + (-1)^{n-r}t^n.
\end{equation}
Comparing (5.18) and (5.19) we get
\begin{multline}
H(t) = (-1)^{n-r}(t^r+\cdots+t^{n-1}) + \\
\sum_{\substack{E\leq n-1\\ e_i\geq 1\text{ for all
$i$}}} \bigg(\sum_{l=0}^n (-1)^{n-l}\binom{n}{l}a^{(l)}_{e_1\dots e_r} 
t^{n-l}\bigg)\frac{p_{e_1}(t)\cdots p_{e_r}(t)}{(1-t)^{E+1}}.
\end{multline}
From the definition of $p_e(t)$ it is straightforward to check by
induction on $e$ that
\begin{equation}
t^{e+1}p_e(1/t) = p_e(t),
\end{equation}
and from the formula (5.10) it is straightforward to check that
\begin{equation}
a^{(n-l)}_{e_1\dots e_r} = (-1)^{n+1+E}a^{(l)}_{e_1\dots
e_r}.
\end{equation}
Equations (5.20), (5.21), and (5.22) imply that
\begin{equation}
t^{n+r-1}H(1/t) = H(t),
\end{equation}
which gives (1.16).

Define a polynomial $g_{e_1\dots e_r}(t)$ by
\begin{align}
g_{e_1\dots e_r}(t) &:= \sum_{l=0}^n (-1)^{n-l}\binom{n}{l}a^{(l)}_{e_1\dots
e_r}t^{n-l} \\
&= (-1)^{n+1+E}\sum_{l=0}^n (-1)^l
\binom{n}{l}a^{(l)}_{e_1\dots e_r}t^l \nonumber
\end{align}
by (5.22).  Then (5.20) gives
\begin{equation}
H(t) = (-1)^{n-r}(t^r + \cdots + t^{n-1}) +
\sum_{\substack{E\leq n-1\\ e_i\geq 1\text{ for all
$i$}}} g_{e_1\dots e_r}(t) \frac{p_{e_1}(t)\cdots
p_{e_r}(t)}{(1-t)^{E+1}}.
\end{equation}
We want to show that $g_{e_1\dots e_r}(t)$ is divisible by
$(1-t)^{E+1}$ and calculate the value at $t=1$ of
$g_{e_1\dots e_r}(t)/(1-t)^{E+1}$.  

Equation (5.10) shows that $a^{(l)}_{e_1\dots e_r}$ is a polynomial in
$l$ of degree $n-1-E$.  It follows that $g_{e_1\dots
e_r}(t)$ is a linear combination of the polynomials
\begin{equation}
\sum_{l=0}^n (-1)^l \binom{n}{l} l^it^l = \bigg(t\frac{d}{dt}\bigg)^i(1-t)^n
\end{equation}
for $i=0,1,\dots,n-1-E$.  These polynomials are clearly all divisible
by $(1-t)^{E+1}$, hence $g_{e_1\dots e_r}(t)$ is also.  

Note that by Taylor's formula, the value at $t=1$ of $g_{e_1\dots
e_r}(t)/(1-t)^{E+1}$ equals the value at $t=1$ of 
\begin{equation}
\frac{(-1)^{E+1}}{(E+1)!}
\bigg(\frac{d}{dt}\bigg)^{E+1} (g_{e_1\dots e_r}(t)).
\end{equation}
Furthermore, the value of this expression at $t=1$ is unchanged if we
replace $(d/dt)$ by $t(d/dt)$.  The polynomial
$(t\frac{d}{dt})^{E+1} g_{e_1\dots e_r}(t)$ is a linear
combination of the polynomials (5.26) for
$i=E+1,\dots,n$.  For $i<n$, the polynomials (5.26) vanish at $t=1$;
for $i=n$, the polynomial (5.26) assumes the value $(-1)^n n!$ at
$t=1$.  Furthermore, equations (5.10) and (5.24) show that when
$(t\frac{d}{dt})^{E+1} g_{e_1\dots e_r}(t)$ is expressed
as a linear combination of the polynomials (5.26), the coefficient of
$(t\frac{d}{dt})^n(1-t)^n$ is
\begin{equation}
\frac{d_1^{e_1}\cdots d_r^{e_r}}{(n-1-E)!\,e_1!\cdots
e_r!}. 
\end{equation}
It follows that the value at $t=1$ of the expression (5.27) is
\begin{equation}
(-1)^{n+1+E}\binom{n}{E+1}
\frac{d_1^{e_1}\cdots d_r^{e_r}}{e_1!\cdots e_r!}.
\end{equation}

It is straightforward to check by induction on $e$ that $p_e(1)=e!$.
From (5.25) and (5.29) we now get
\begin{equation}
H(1) = (-1)^{n-r}(n-r) + \sum_{\substack{E\leq n-1\\
e_i\geq 1\text{ for all $i$}}} (-1)^{n+1+E}
\binom{n}{E+1} d_1^{e_1}\cdots d_r^{e_r},
\end{equation}
which is (1.15).

\section{The case $r\geq n$}

Let $C^\bullet(f_1,\dots,f_r)$ be the (cohomological) Koszul complex
on $K[x]$ defined by $f_1,\dots,f_r$.  Consider the grading on $K[x]$
defined by total degree in $x_1,\dots,x_n$ and let $K[x]^{(i)}$ denote
the space of homogeneous polynomials of degree $i$.  This induces a
grading on $C^\bullet(f_1,\dots,f_r)$ by defining
\[ C^k(f_1,\dots,f_r)^{(i)} = \bigoplus_{1\leq j_1<\dots<j_k\leq r}
K[x]^{(i+d_{j_1}+\cdots+d_{j_k})}, \]
i.~e., the grading is determined by requiring that 
$C^0(f_1,\dots,f_r)^{(i)} = K[x]^{(i)}$ and that the boundary maps
are graded homomorphisms.  The following lemma is probably well known,
but we do not know a reference for it.

\begin{lemma}
Suppose the ideal $(f_1,\dots,f_r)$ has depth $n$, i.~e.,
$f_1,\dots,f_r$ have no common zero in ${\bf P}^{n-1}$.  Then
\begin{equation}
H^k(C^\bullet(f_1,\dots,f_r)^{(i)}) = 0\quad\text{for $i>-n$ and all
$k$}
\end{equation}
and
\begin{equation}
\dim_KH^k(C^\bullet(f_1,\dots,f_r)^{(-n)}) = \begin{cases} 1 &
\text{if $k=n$,} \\ 0 & \text{otherwise.} \end{cases}
\end{equation}
\end{lemma}

{\it Proof}.  We prove the result for the ideal
$(x_1,\dots,x_n,f_1,\dots,f_r)$ and then explain how to inductively
remove $x_1,\dots,x_n$.  It is well known that
$H^k(C^\bullet(x_1,\dots,f_r))$ is isomorphic to the cohomology of the
Koszul complex on $K[x]/(x_1,\dots,x_n)(\cong K)$ defined by
$f_1,\dots,f_r$.  Let $\overline{C}{}^\bullet$ denote this latter
Koszul complex.  In terms of the gradings, we have more precisely
\begin{equation}
H^{k+n}(C^\bullet(x_1,\dots,f_r)^{(i)})\cong
H^k((\overline{C}{}^\bullet)^{(i+n)}), 
\end{equation}
where 
\[ (\overline{C}{}^0)^{(i)} = \begin{cases} K & \text{if
$i=0$,} \\ 0 & \text{otherwise.} \end{cases} \]
One has trivially
\[ H^k((\overline{C}{}^\bullet)^{(i)}) = 0\quad\text{for $i>0$ and all
$k$} \]
and
\[ \dim_K H^k((\overline{C}{}^\bullet)^{(0)}) = \begin{cases} 1 &
\text{if $k=0$,} \\ 0 & \text{otherwise,} \end{cases} \]
so the assertions of the lemma for $(x_1,\dots,f_r)$ follow from (6.4).

For notational convenience, put 
\[ C_l^\bullet = C^\bullet(x_1,\dots,x_l,f_1,\dots,f_r). \]
Suppose inductively the assertions of the lemma are true for some
$C_l^\bullet$, where $1\leq l\leq n$.  We prove them for
$C_{l-1}^\bullet$.  There is a well known short exact sequence of
graded Koszul complexes (see \cite[Theorem 16.4]{MA})
\[ 0\to (C_{l-1}^\bullet)^{(i+1)}[-1]\to (C_l^\bullet)^{(i)}\to
(C_{l-1}^\bullet)^{(i)}\to 0, \]
which gives rise to the exact cohomology sequence
\begin{equation}
\ldots\to H^k((C_l^\bullet)^{(i)})\to H^k((C_{l-1}^\bullet)^{(i)})\to
H^k((C_{l-1}^\bullet)^{(i+1)})\to
H^{k+1}((C_l^\bullet)^{(i)})\to\dots.
\end{equation}
By the induction hypothesis,
\[ H^k((C_l^\bullet)^{(i)}) = 0\quad\text{for $i>-n$ and all $k$,} \]
so we get isomorphisms
\begin{equation}
H^k((C_{l-1}^\bullet)^{(i)})\cong H^k((C_{l-1}^\bullet)^{(i+1)})
\quad\text{for $i>-n$ and all $k$.}
\end{equation}
The graded cohomology groups $H^k(C_{l-1}^\bullet)$ are annihilated by the
ideal $(f_1,\dots,f_r)$, so our hypothesis implies that $x_1,\dots,x_n$
are contained in the radical of the annihilator of this graded
module.  It follows that these cohomology groups are
finite-dimensional, hence
\[ H^k((C_{l-1}^\bullet)^{(i)})=0\quad\text{for $i$ sufficiently large
and all $k$.} \]
Using (6.6) and descending induction on $i$, we get
\begin{equation}
H^k((C_{l-1}^\bullet)^{(i)})=0\quad\text{for $i>-n$ and all $k$.}
\end{equation}
Taking $i=-n$ in (6.5) and using (6.7) now gives
\[ H^k((C_l^\bullet)^{(-n)})\cong H^k((C_{l-1}^\bullet)^{(-n)})
\quad\text{for all $k$,} \]
thus the assertions of the lemma hold for $C_{l-1}^\bullet$.

The main result of this section is the following.
\begin{proposition}
Suppose that $f_1,\dots,f_r$ have no common zero in ${\bf P}^{n-1}$.
Then
\begin{equation}
H^k(\Omega^\bullet_{K[x,y]/K})^{(0,p)} = 0\quad\text{for $k\neq 2n$
and all $p$,}
\end{equation}
and
\begin{equation}
\dim_K H^{2n}(\Omega^\bullet_{K[x,y]/K})^{(0,p)} = \begin{cases} 1 &
\text{if $p=n$} \\ 0 & \text{otherwise.} \end{cases}
\end{equation}
\end{proposition}

{\it Proof}.  Regard $(\Omega^\bullet_{K[x,y]/K})^{(0,\bullet)}$ as
the total complex associated to the double complex whose vertical map
$\partial_v$ is the wedge product with $\sum_{j=1}^r f_jdy_j$ and
whose horizontal map $\partial_u$ is the wedge product with
$\sum_{j=1}^r y_jdf_j$.  The $l$-th column of this double complex is
the direct sum over $y_1^{b_1}\cdots
y_r^{b_r}\,dx_{i_1}\cdots dx_{i_l}$ of complexes whose component in row
$m$ is 
\[ \bigoplus_{1\leq j_1<\dots<j_m\leq r}
K[x]^{(d_{j_1}+\cdots+d_{j_m}-l+\sum_{j=1}^r b_jd_j)}y_1^{b_1}\cdots
y_r^{b_r}\,dx_{i_1}\cdots dx_{i_l}dy_{j_1}\cdots dy_{j_m}. \]
This complex is clearly isomorphic to
$C^\bullet(f_1,\dots,f_r)^{(-l+\sum_{j=1}^r b_jd_j)}$.  By Lemma 6.1,
all the vertical cohomology vanishes unless $b_j=0$ for all $j$ and
$l=n$.  In this latter case, the lemma implies that all vertical
cohomology vanishes except in row $n$ of column $n$, where it is
one-dimensional.  The proposition now follows by computing the
cohomology of the total complex as the horizontal cohomology of the
vertical cohomology.

We give an explicit basis for
$H^{2n}(\Omega^\bullet_{K[x,y]/K})^{(0,n)}$ when $d_1\cdots d_r\neq 0$
in $K$.  Let $\xi_n$ be the $2n$-form given by~(4.23).
\begin{lemma}
$\partial(\xi_n) = 0$
\end{lemma}

{\it Proof}.  It is easily seen that $\partial(\xi_n) =
(f_1dy_1+\cdots+f_rdy_r)\wedge\xi_n$ and that the coefficient of
$dx_1\wedge\cdots\wedge dx_n\wedge dy_{i_1}\wedge\cdots\wedge
dy_{i_{n+1}}$ in $\partial(\xi_n)$ is, up to sign,
\[ \bigg(\prod_{i\not\in\{i_1\dots,i_{n+1}\}}d_i\bigg)\det\begin{bmatrix}
d_{i_1}f_{i_1} & \frac{\partial f_{i_1}}{\partial x_1} &
\cdots & \frac{\partial f_{i_1}}{\partial x_n} \\
d_{i_2}f_{i_2} & \frac{\partial f_{i_2}}{\partial x_1} &
\cdots & \frac{\partial f_{i_2}}{\partial x_n} \\
\hdotsfor{4} \\
d_{i_{n+1}}f_{i_{n+1}} & \frac{\partial f_{i_{n+1}}}{\partial x_1} &
\cdots & \frac{\partial f_{i_{n+1}}}{\partial x_n}
\end{bmatrix}. \]
By the Euler relation, the first column is a $K[x]$-linear combination
of the other columns, hence this determinant is zero. 

\begin{proposition}
If $d_1\cdots d_r\neq 0$, then $[\xi_n]$ is a basis for
$H^{2n}(\Omega^\bullet_{K[x,y]/K})^{(0,n)}$. 
\end{proposition}

{\it Proof}.  Using (4.17), (4.12), (6.9), and (6.10), one proves analogues of
(4.19) and (4.20) by induction on $k$: 
\begin{gather}
H^{2k}(\widetilde{\Omega}^\bullet_{K[x,y]/K})^{(0,p)} = 0 \quad \text{for
$0\leq k<n$ and all $p$} \\
\intertext{and}
H^{2k-1}(\widetilde{\Omega}^\bullet_{K[x,y]/K})^{(0,p)} \cong \left.
\begin{cases} K & \text{if $p=k$} \\ 0 & \text{otherwise} \end{cases}
\right\} \quad\text{for $1\leq k\leq n$.}
\end{gather}
As in the proof of Proposition 4.24, one has 
\[ \theta(\xi_n) = (-1)^{n(n-1)/2}(d_1\cdots d_r)\eta_n +
dF\wedge\theta(\tau_n). \]
Since $[\eta_n]$ is a basis for
$H^{2n-1}(\widetilde{\Omega}^\bullet_{K[x,y]/K})^{(0,n)}$ and
$d_1\cdots d_r\neq 0$, $[\xi_n]$
is not trivial in $H^{2n}(\Omega^\bullet_{K[x,y]/K})^{(0,n)}$.  By
(6.10), $[\xi_n]$ must be a basis for
$H^{2n}(\Omega^\bullet_{K[x,y]/K})^{(0,n)}$. 

{\it Remark}.  Note that when $r=n$,
\[ \xi_n = \frac{\partial(f_1,\dots,f_n)}{\partial(x_1,\dots,x_n)}\, 
dx_1\wedge\cdots\wedge dx_n\wedge dy_1\wedge\cdots\wedge dy_n, \]
where $\partial(f_1,\dots,f_n)/\partial(x_1,\dots,x_n)$ denotes the
Jacobian determinant.  In this case, the nontriviality of $[\xi_n]$
in $H^{2n}(\Omega^\bullet_{K[x,y]/K})^{(0,n)}$ is equivalent to the
assertion that
\[ \frac{\partial(f_1,\dots,f_n)}{\partial(x_1,\dots,x_n)}\not\in
(f_1,\dots,f_n). \]
The earliest reference for this fact of which we are aware is \cite[Corollary
4.7]{SS}.  If $K$ has characteristic zero, it can be proved by a residue
argument (see \cite[Chapter~5, Section 2]{GH} or \cite[Theorem 12.6(ii)]{PS}). 
More generally, it follows from Proposition~6.12 that a basis for
$H^n(C^\bullet(f_1,\dots,f_r)^{(-n)})$ (see (6.3)) is the cohomology
class of the element of $C^n(f_1,\dots,f_r)^{(-n)}$ whose component
corresponding to an $n$-tuple $1\leq j_1<\cdots<j_n\leq r$ is the
Jacobian determinant
$\partial(f_{j_1},\dots,f_{j_n})/\partial(x_1,\dots,x_n)$.

\end{document}